\documentclass[a4paper,10pt,twocolumn]{article}
\usepackage[utf8]{inputenc}

\usepackage[switch]{lineno}

\usepackage{hyperref}
\usepackage{bm}
\usepackage{authblk}
\usepackage{siunitx}
\usepackage{graphicx}

\usepackage[a4paper,left=1.5cm,right=1.5cm,top=2cm,bottom=2cm]{geometry}

\usepackage{nomencl}
\usepackage{ifthen}
\renewcommand{\nomgroup}[1]{%
\ifthenelse{\equal{#1}{C}}{\item[\textbf{Constants}]}{%
\ifthenelse{\equal{#1}{G}}{\item[\textbf{Greek Letters}]}{%
\ifthenelse{\equal{#1}{S}}{\item[\textbf{Subscripts}]}{}}}
}
\makenomenclature

\nomenclature{H}{Blade Height [m]}
\nomenclature{D}{Rotor Diameter [m]}
\nomenclature{M}{Moment Exerted On The Rotor Shaft [N.m]}
\nomenclature{$\Omega$}{Rotational Velocity  [m]}
\nomenclature{$C_P$}{Power Coefficient}
\nomenclature{c}{Airfoil Chord [m]}
\nomenclature{U}{Streamwise Component Of Velocity Vector [m/s]}
\nomenclature{V}{Transversal Component Of Velocity Vector [m/s]}
\nomenclature{Re}{Reynolds Number }
\nomenclature{t}{Time [sec]}
\nomenclature{$\rho$}{Density $[kg/m^3]$}
\nomenclature{$\mu$}{Dynamic Viscosity [Pa.s]}
\nomenclature{$\nu$}{Kinematic Viscosity $[m^2/s]$}
\nomenclature{$\Phi_U$}{Left Eigenvector}
\nomenclature{$\Phi_V$}{Right Eigenvector}
\nomenclature{$\Sigma$}{Eigenvectors}

\usepackage[authoryear]{natbib}

\begin{document}
\let\WriteBookmarks\relax
\def\floatpagepagefraction{1}
\def\textpagefraction{.001}

\title{Assessment of URANS and LES Methods in Predicting Wake Shed Behind a Vertical Axis Wind Turbine}
\author[1]{Armin~Sheidani\footnote{armin.sheidani@polimi.it}}
\author[2]{Sajad~Salavatidezfouli\footnote{ssalavat@sissa.it}}
\author[2]{Giovanni~Stabile\footnote{gstabile@sissa.it}}
\author[2]{Gianluigi~Rozza\footnote{grozza@sissa.it}}

\affil[1]{Department of Mechanical Engineering, Politecnico di Milano, Milan, Italy}
\affil[2]{Mathematics Area, mathLab, SISSA, via Bonomea 265, I-34136 Trieste, Italy}

\date{} 

\twocolumn[
  \begin{@twocolumnfalse}
    \maketitle
	\begin{abstract}
		In order to shed light on the Vertical-Axis Wind Turbines (VAWT) wake characteristics, in this paper we present high-fidelity CFD simulations of the flow around an exemplary H-shaped VAWT turbine, and we propose to apply Proper Orthogonal Decomposition (POD) to the computed flow field in the near wake of the rotor. The turbine under consideration was widely studied in previous experimental and computational investigations. In the first part of the study, multiple Reynolds-Averaged Navier-Stokes (RANS) simulations were performed at the Tip Speed Ratio (TSR) of peak power coefficient, to select the most accurate turbulence model with respect to available data. In the following step, further RANS numerical simulations were performed at different TSRs to compare the power coefficient against experimental data. Then, Large Eddy Simulation (LES) was applied for multiple TSR conditions. The spatial and temporal POD modes along with modal energy for the RANS and LES results were extracted, and the performance of the turbulence models was assessed. Also, an interpretation of the POD modes with respect to the flow structures was given to highlight the most significant time and length scales of the predictions considering the different dynamical levels of approximations of the computational models.
		
		\vspace{0.5cm}
		\textbf{Highlights:}
		\begin{itemize}
		  \item LES method computed the wake region characteristics more accurately compared to RANS.
		
		  \item RANS method captures the wake region length and width close to that of the LES results.
		
		  \item After the third POD mode, the RANS method fails in accurately simulating flow structures.
		
		\end{itemize}
		
		\end{abstract}
  \end{@twocolumnfalse}
]
\maketitle
\printnomenclature

\section{Introduction}

Lift-driven Vertical-Axis Wind Turbines (VAWT) exhibit interesting advantages with respect to the widespread horizontal-axis counterpart in terms of installation and operation \citep{Tjiu_2015a, Tjiu_2015b}, even though they are characterized by complex unsteady aerodynamics which eventually penalizes the energy conversion efficiency of the rotor \citep{elcheikh2018performance, Chen_2022}. Typically, VAWT operation is inherently associated with wide oscillations in Reynolds number and angle of attack during the revolution which makes the blades undergo dynamic-stall and even deep-stall \citep{Sheidani_2022, Tirandaz_2021}. As a result of such oscillations, vortices detaching from the blades are shed into the wake, making it highly unsteady and composed of multiple time/length scales. These features make VAWT wakes highly different from those of horizontal-axis rotors \citep{Hamlaoui_2022}, with important implications on the mixing rate and evolution of the flow downstream of the turbine. 

Owing to the unique characteristics of lift-driven VAWT such as low sensitivity to the wind direction, relatively higher power coefficient compared to that of drag-type VAWT \citep{Franchina_2020} and, more importantly, smaller dimensions, VAWTs are growing in popularity to be utilized on a large scale ranging from urban to off-shore areas \citep{stathopoulos2018urban, Franchina_2019}, following with consequent studies on the compactness and optimizing tandem configurations \citep{Silva_2021}. Therefore, the presence of the turbines in large numbers and denser configurations, heightens the strong need to gain a deeper insight into the wakes shed past a VAWT \citep{hohman2020effect, vergaerde2020experimental, posa2020dependence, posa2020influence}. 
 
 As a whole, a deep inspection into the flow characteristics may be obtained by means of high-fidelity experimental and numerical methods such as Particle Image Velocimetry (PIV) \citep{Lee_2016, Sim_o_Ferreira_2008, Herr_ez_2018, Ferreira_2020} and Scale Resolving Simulation (SRS). While resolving the full spectrum of turbulence by means of Direct Numerical Simulation (DNS) remains out of reach due to excessive computational cost, SRS methods including Large Eddy Simulation (LES) \citep{Safari_2018, He_2021, Tang_2022} and Detached Eddy Simulation (DES) \citep{Larin_2016, Boudreau_2017, Qian_2019} and Scale-Adaptive Simulation (SAS) \citep{Toja_Silva_2018, Rezaeiha_2019b} aim to resolve only the larger scales of the turbulence with modelling the smaller ones. However, these techniques impose a great deal of computational cost to the system as well and yet they are prone to the existence of noise which requires employment of different methods to eliminate the unfavorable indeterministic signals \citep{Scherl_2020}. As a result, this study, as one of the primary objectives, is intended to investigate to what extent more economical turbulence models like Reynolds Average Navier Stokes (RANS) can be implemented to study the phenomenon of wake shed behind a VAWT. For this purpose, Proper Orthogonal Decomposition (POD) will be utilized to extract the hidden structures of both turbulence simulation techniques to cast light deeply on this matter.
 
The employment of POD to extract different modes of flow dynamics was first proposed by Lumely \citep{lumley1967structure}. This idea was then nurtured and developed into more robust methods in different studies namely, Sirovich \citep{Sirovich_1987} who put forward the concept of snapshot POD. Consequently, in recent years, the advances in Reduced Order Modelling (ROM) techniques have led to a great deal of studies intended to enlighten the hidden flow characteristics by employing methods like POD \citep{aroma_book,lassila2014model, _tefan_2017, Sinha_2018, Wang_2018, Chouak_2018, benner2021system, benner2020modela, benner2020modelb}.

There have been several attempts to investigate wake characteristics shed behind a VAWT. In this regard, Abkar and Dabiri \citep{Abkar_2017} employed the LES method to study the wake of a Darrieus VAWT. It was shown that the wakes of the VAWT demonstrate self-similarity especially, in the central part of the wake. Tescione et al. \citep{Tescione_2014} performed an experimental study by implementing PIV technique, to study the wake characteristics of an H-shaped VAWT at a close distance from the rotor. One of the main observations in this study was that the vortical structures can be divided into two major parts mainly at a distance being equal to three times the rotor radius. Kuang et al. \citep{Kuang_2022} employed Improved Delayed Detached Eddy Simulation (IDDES) to investigate the effect of the wake of an upstream floating VAWT on a downstream one. It was reported that the pitch motion leads to a lower amount of velocity deficit in the core region of the wake area. Shamsoddin and Porté-Agel \citep{Shamsoddin_2020} studied the effect of aspect ratio on the wake shed behind a VAWT and it was shown that while the wake keeps the primary aspect ratio in the near wake region, it tends to an elliptical shape as it moves away from the rotor. 

Owing to the distinctive wake characteristics of VAWTs there have been several studies aiming at deriving an analytical model for the description of the VAWT wake. In this regard, Peng et al. \citep{Peng_2021} reviewed different analytical models validated by a wide range of high-fidelity methods including field test, PIV and CFD methods for this purpose. It should be noted that the change in wake characteristics behind a VAWT due to the blade geometrical alterations has been recently studied \citep{Asim_2021}. High-fidelity CFD models often contain too much information at the cost of severe computationals. However, despite the lower accuracy associated with RANS methods compared to LES, as found in the literature \citep{HijaziStabileMolaRozza2020b, Lorenzi_2016, Nakazawa_2019, Liu_2014}, the employment of POD on RANS results yet yields interpretable outcome which may be utilized to extract valuable flow characteristics. Chen et al. \citep{Chen_2020} investigated the effect of gurney flaps on the performance improvement of an H-shaped Darrieus VAWT by means of RANS turbulence models. Following this, a POD analysis on the turbine wake was implemented and it was shown that the use of gurney flaps leads to curbing the vortex strength in the wake of the turbine after the fifth POD mode.  

 Considering the studies conducted so far on Reduced Order Modelling of flow characteristics past a wind turbine, it reveals that the modal analysis of wake shed behind a Darrieus VAWT rotor, by highly accurate methods like LES, in order to unfold the major flow structures has remained unnoticed. More importantly, with regards to the work of Asim and Islam \citep{Asim_2021} and Chen et al. \citep{Chen_2020} where the POD analysis has been separately performed on LES and RANS results of Savenious and Darrieus wind turbines simulation respectively, this study is intended to assess to what extent the RANS methods can be trusted to perform the wake study or extract the small flow structures with respect to LES modes considered as the basis for the comparison. This point is of great importance owing to the growing popularity of VAWTs in more compact configurations leading to the necessity of performing an in-depth study on the flow characteristics past a VAWT. Therefore, in this study, in order to perform an analysis on the VAWT wake, POD method will be employed to unfold the physics of flow behind a VAWT. In the beginning, the numerical methods employed in this study will be discussed. Following this, the results will be validated with respect to the experimental study of Battisti et al. \citep{Battisti_2018} and the POD method will be briefly discussed. Finally, a discussion on the LES and RANS results and the summary will be presented.

\section{NUMERICAL METHODS}

In order to perform the numerical simulation of the problem described in the previous section, the Unsteady Reynolds Average Navier Stokes (URANS) and the filtered form of Navier Stokes for the LES equations are numerically solved using ANSYS FLUENT 2020 software. 
The URANS formulation is adopted for the continuity and momentum equations of the incompressible fluid which is as follows:

\begin{equation}
\label{eq1}
\frac{\partial \overline{U_{i}}}{\partial x_i}=0\\
\end{equation}

\begin{equation}
\label{eq2}
\frac{\partial \rho \overline{U_i}}{\partial t}+\frac{\partial \rho \overline{U_i \bar{U}_j}}{\partial x_j}=-\frac{\partial P}{\partial x_i}+\frac{\partial}{\partial x_j}\left[\mu\left(\frac{\partial \overline{U_i}}{\partial x_j}+\frac{\partial \overline{U_j}}{\partial x_i}\right)-\rho \overline{u_{i}^{\prime} u_{j}^{\prime}}\right]\\
\end{equation}

where $\overline{U_{i}}$ and $\rho \overline{u_{i}^{\prime} u_{j}^{\prime}}$ are the averaged velocity and Reynolds stress tensor, respectively. The latter can be written as:

\begin{equation}
\label{eq3}
-\rho \overline{\boldsymbol{u}_{i}^{\prime} \boldsymbol{u}_{\boldsymbol{j}}^{\prime}}=2 \mu_t \overline{S_{ij}}-\frac{2}{3} k \rho \overline{\delta_{ij}}\\
\end{equation}

where the mean strain rate tensor is:

\begin{equation}
\label{eq4}
\overline{S_{ij}}=\frac{1}{2}\left(\frac{\partial \bar{U}_i}{\partial x_j}+\frac{\partial \bar{U}_j}{\partial x_i}\right)\\
\end{equation}

The turbulent viscosity term, $\mu_t$, can be addressed with common one- or two- equation eddy viscosity models to provide turbulence closure.  
As for the LES model, filtered Navier-Stokes is utilized:

\begin{equation}
\label{eq5}
\frac{\partial \rho \widetilde{U_i}}{\partial t}+\frac{\partial \rho \widetilde{U_i} \widetilde{U_j}}{\partial x_j}=-\frac{\partial P}{\partial x_i}+\frac{\partial}{\partial x_j}\left[\mu\left(\frac{\partial \widetilde{U_i}}{\partial x_j}+\frac{\partial \widetilde{U_j}}{\partial x_i}\right)\right]-\frac{\partial \tau_{\mathrm{ij}}}{\partial x_j}\\
\end{equation}

where the last term, $\tau_{\mathrm{ij}}$, is the sub-grid scale (SGS) stress and defined as follows:

\begin{equation}
\label{eq6}
\tau_{i j}=\rho \widetilde{u_i u_j}-\rho \widetilde{u_i} \widetilde{u_j}
\end{equation}

There exists couple of closure models to address SGS in terms of the local resolved flow. Selection of the appropriate turbulence closure model for RANS and LES will be discussed further on. In the following section the structure of the numerical method including the grid generation, turbulence modelling techniques, POD method, independency studies and validation results will be presented.

\subsection{Problem Description}

\begin{figure*}
	\centering
	\includegraphics[width=\textwidth]{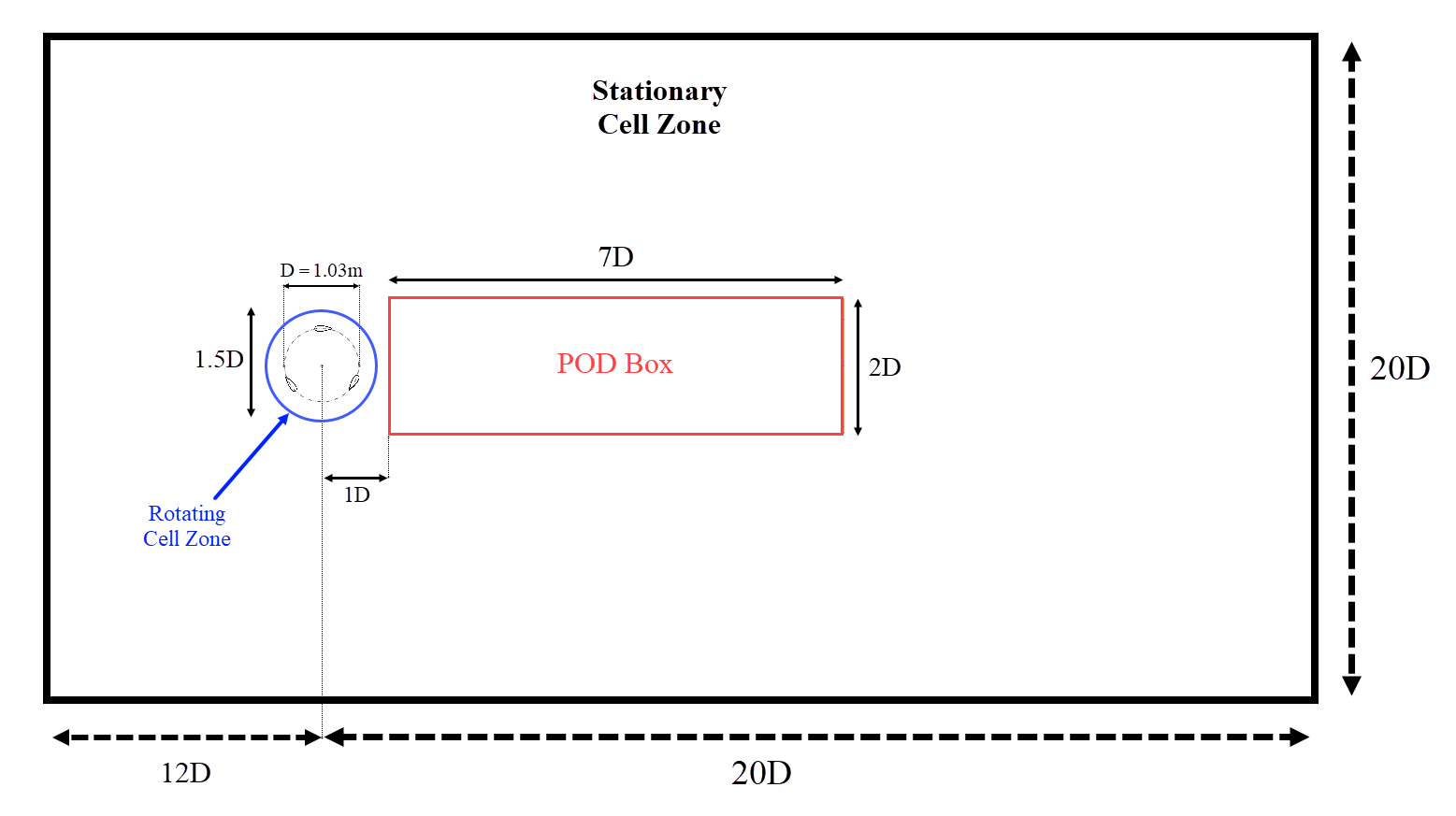}
	\caption{Dimensions of the cell zones of the computational domain}
	\label{FIG1}
\end{figure*}

This study is aimed to perform the RANS and LES numerical simulations on a VAWT based on the experimental work of Battisti et al. \citep{Battisti_2018} at different tip speed ratios (TSR). The rotor diameter of the VAWT is equal to 1.03 (m) and it consists of three NACA0021 airfoils with a chord value of 0.086 m. The simulations will be implemented at the rated TSR being equal to 2.4 and the off-design points of operation corresponding to TSR3.3 and TSR1.5. It should be noted that as the turbine rotates at 400(rpm) the TSRs 1.5, 2.4 and 3.3 correspond to inlet velocity of $14.38 (m\slash s)$, $9.1 (m\slash s)$ and $6.5 (m\slash s)$ respectively. Also, this study is intended to perform POD analysis and investigate the wakes shed behind the rotor as shown in the box presented in Fig. \ref{FIG1} where the dimensions used in the computational domains have been designated. The inlet and outlet boundaries are respectively on the left and right sides of the domain and the upper and lower boundaries are considered to be of pressure outlet type.

\subsection{Grid Generation} 
The necessity of the full-scale simulation of the turbine including blades and support structure (support arms and central hub) is dependent on the type of the wind turbine (\cite{sotoudeh2019field,mousavi2020mathematical}). Siddiqui et al. \citep{siddiqui2015quantification} concluded that for this specific type of VAWT, i.e. H-Darrieus wind turbine, effects of the upper and lower tips are negligible. Hence, employment of a 3D model composed of a portion of the blades is sufficiently accurate. To this end, 3D meshes were created by extrusion of a 2D mesh into the normal direction by means of one chord length, as recommended by \citep{xu2017delayed,li20132}, while considering the periodic boundary condition for two sides of the extrusion.

In order to perform a grid independency study, three sets of meshes with $1\times 10^6$, $2.1\times 10^6$ and $3.8\times 10^6$ hexahedral elements were created, respectively, for the RANS simulations. There observed less than a 1\% difference between the results of two later meshes, and hence, the second mesh was selected to perform the post-processing.

\begin{figure*}
\centering
\includegraphics[scale=.21]{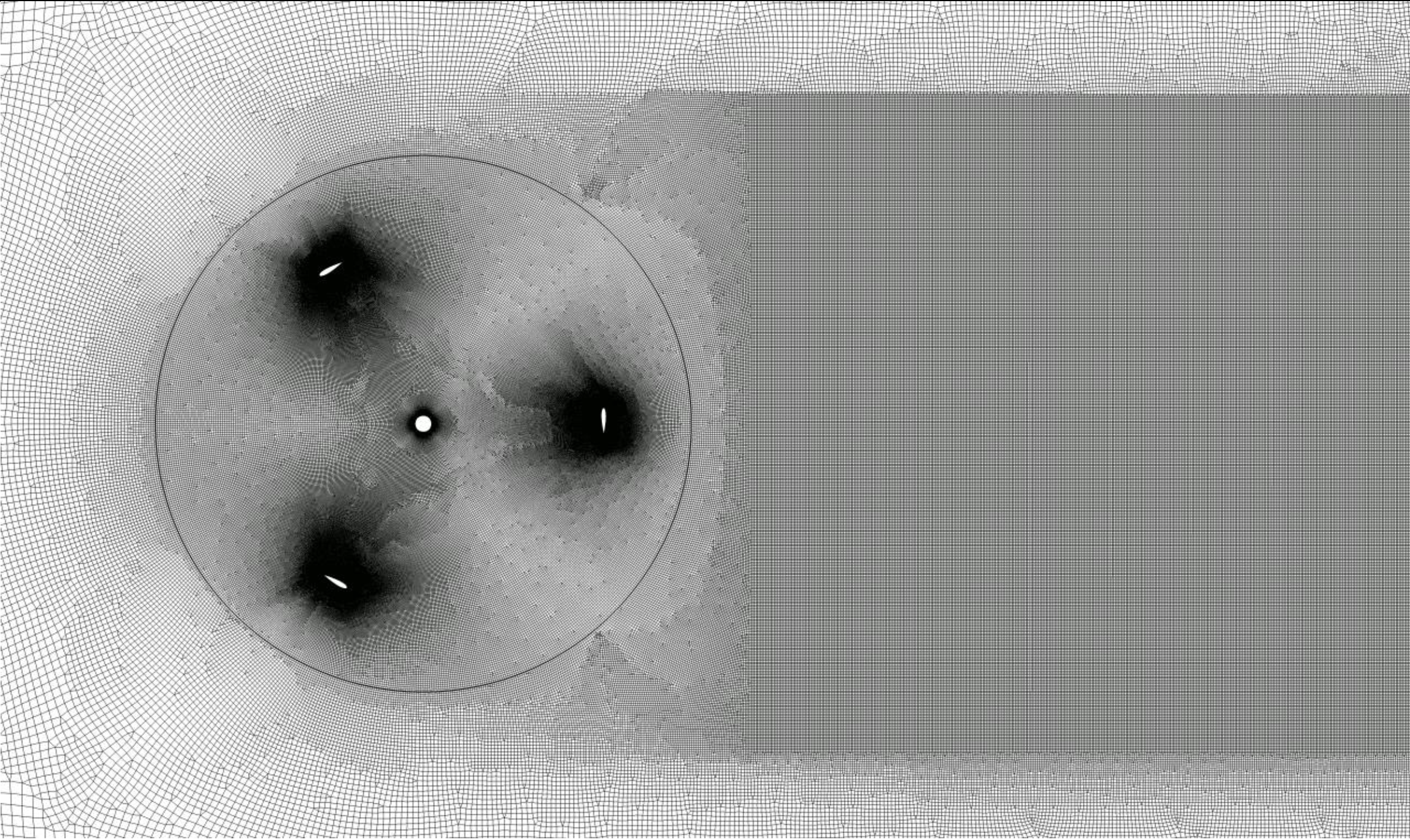}
\includegraphics[scale=.21]{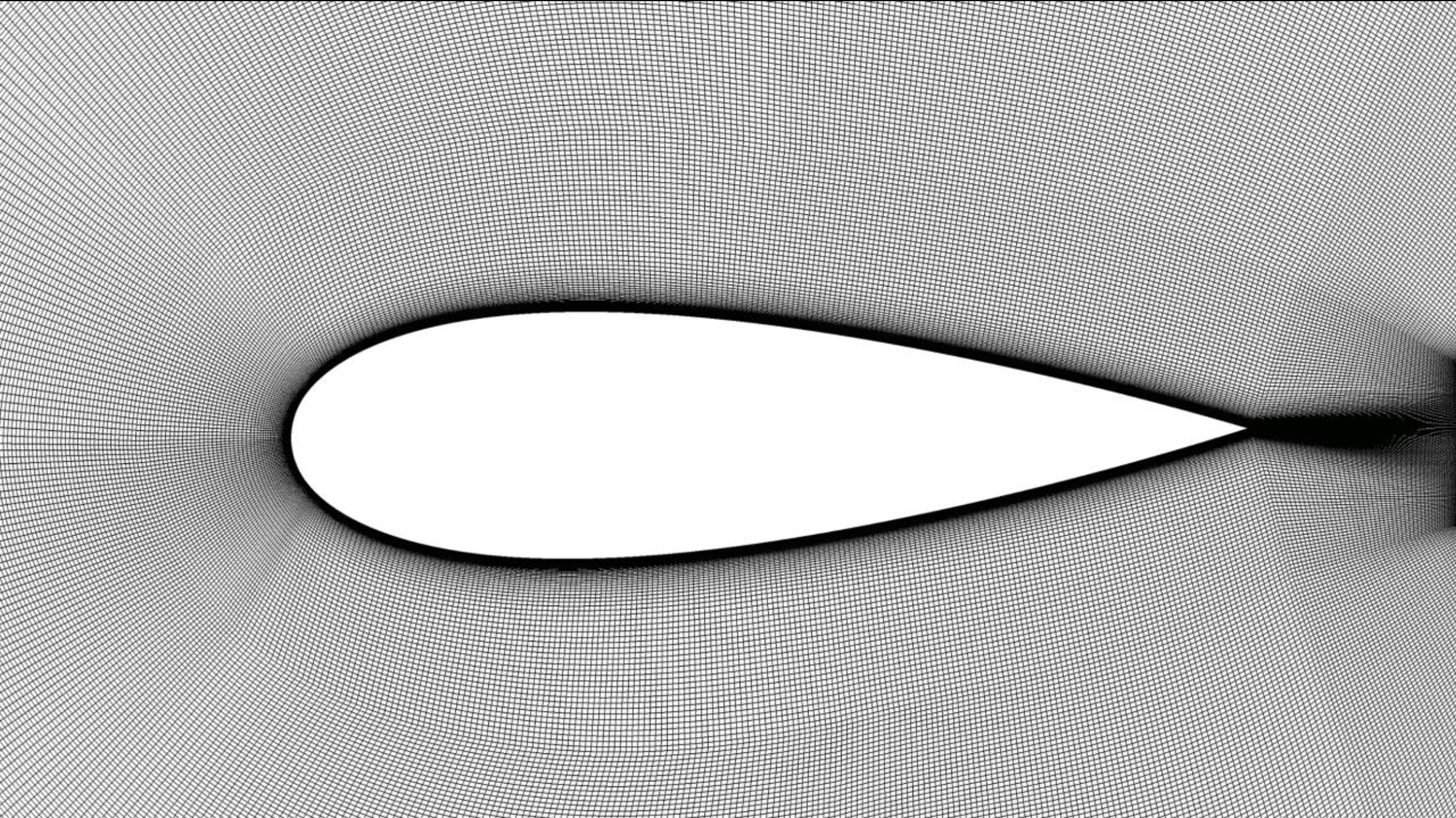}
\includegraphics[scale=.1]{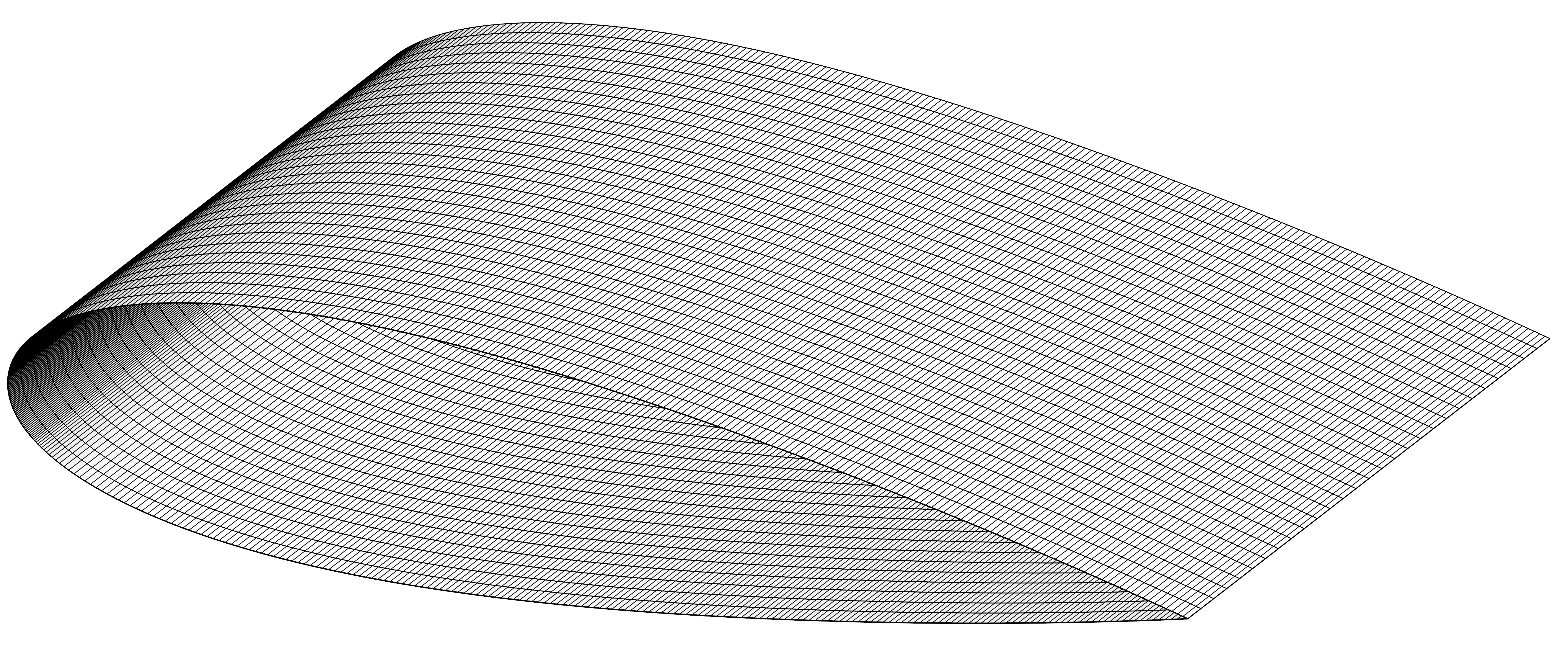}
\caption{ Overview of the LES mesh (a) around the turbine and consequent cell zones b) around the blade c) on the surface of the extruded blade.}
\label{FIG2}
\end{figure*}

\begin{figure*}
\centering
\includegraphics[scale=.2]{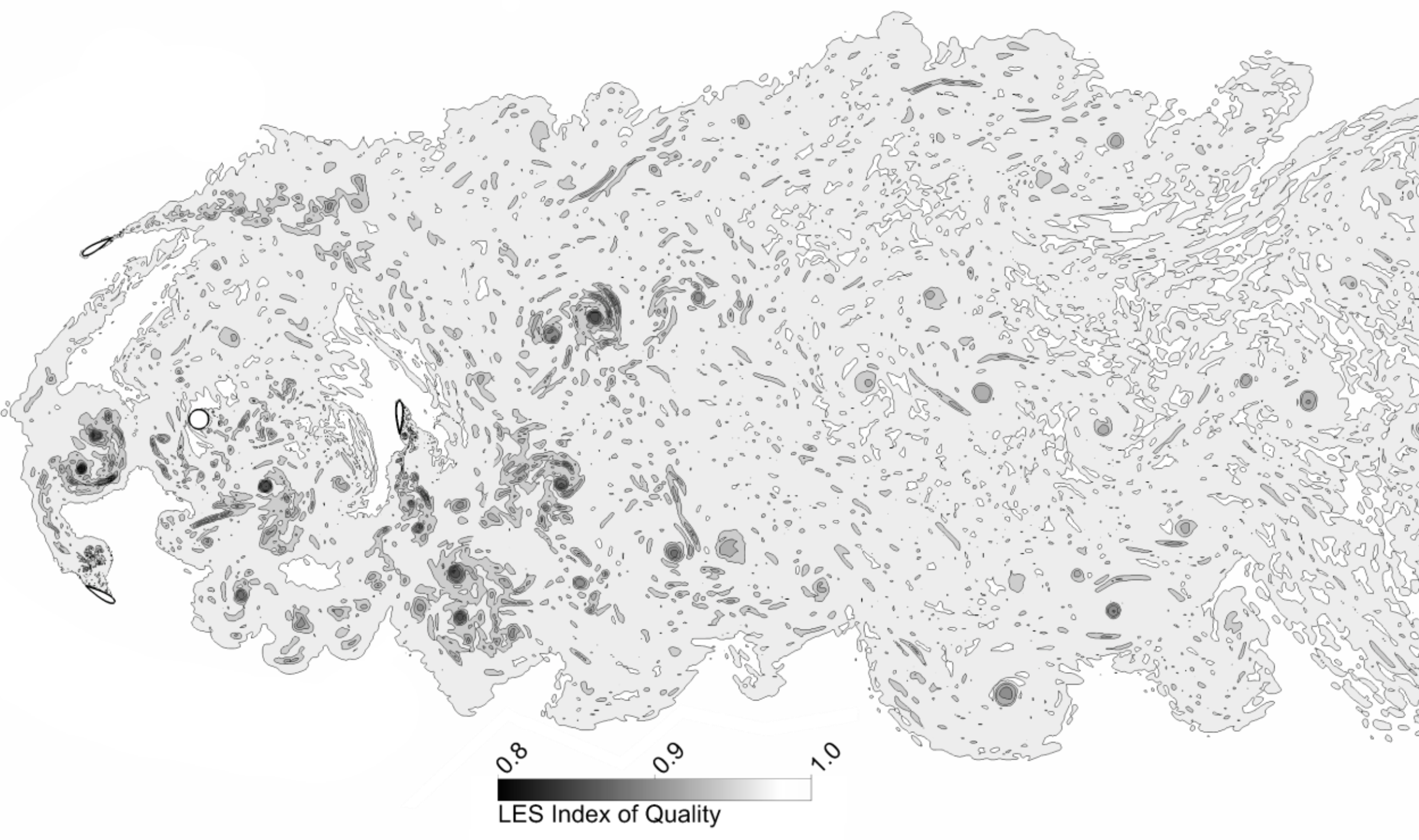}
\caption{Contour of \textit{index quality} for the LES mesh}
\label{FIG3}
\end{figure*}

As for the LES, the model consisted of $1.1\times 10^7$ hexahedral elements with a total of 28 layers in the spanwise direction. The minimum and average orthogonal qualities of the grid were 0.3 and 0.9, respectively, which is regarded as a good quality for LES simulations \citep{yagmur2021numerical}. The minimum resolution percentage based on the Index Quality criterion proposed by Celik et al \citep{Celik_2005} was found to be 81\% (Fig. \ref{FIG3}). According to \citep{pope2000turbulent, Celik_2005} the value above 80\% corresponds to an adequate resolved portion of the turbulent kinetic energy for the LES simulation indicating an appropriate mesh density. Finally, the first element height for wall boundaries was set to $5\times 10^{-6}$ corresponding to $y^{+}\sim1$. The schematic presentation of the LES grid can be seen in Fig. \ref{FIG2}.

\subsection{Turbulence Modelling Techniques}

As the blades during one period of rotation undergo a wide range of Reynolds number from very low to high, the suitable RANS model needs to be found. To this end several fully turbulent and transient models including $\gamma-Re_{\theta}$, $k_{t}-k_{l}-\omega$ and $k-\omega$ SST were performed at the rated TSR and it was found that the latter i.e. $k-\omega$ SST yields the most accurate results. It should be noted that this is matching the finding of Rezaeiha et al. \citep{Rezaeiha_2019}. In addition, all the spatial and temporal terms in Equation \ref{eq2} were discretized by second order method and SIMPLE algorithm was employed to solve the set of equations. Furthermore, regarding the LES equations, Wall-Adapting Local Eddy-Viscosity (WALE) Model \citep{nicoud1999subgrid} was implemented to capture the Sub-Grid Scale (SGS) structures as suggested by Ma et al. \citep{Ma_2009} to be the most suitable modelling technique for the LES computation of the flow physics of wind turbines. In addition, the convection term in Equation \ref{eq5} was discretized by Central Differencing method as suggested by Li et al. \citep{Li_2013} where, similar to the current study, a 2.5D LES of a VAWT had been implemented. Also, the time step size was set to be $10^{-4} sec$ as it corresponds to $0.24^{\circ}$ of rotation leading to a CFL number less than unity.

\subsection{Proper Orthogonal Decomposition Method}

Proper Orthogonal Decomposition is employed to extract the coherent structures both in spatial and temporal domain \citep{Li_2013}. It should be noted that the POD method implemented in this study is based on the work of Taira et al. \citep{Taira_2017}. In order for POD to be implemented, the Singular Value Decomposition (SVD) method was employed owing to its robustness against roundoff errors \citep{Lengani_2014}. In order to perform SVD, the 'matrix of snapshots' denoted as 'X' needs to be formed. The rows of X include the value of the variable of interest corresponding to its position in the computational domain and the columns correspond to the time instance at which the data has been extracted. Having formed the X matrix, SVD is applied as follows in Equation 3-7:

\begin{equation}
\boldsymbol{X}=\Phi_{\mathrm{U}} \Sigma \Phi_{\mathrm{V}}^T\\
\end{equation}

Where $\Phi_{\mathrm{U}}$ and $\Phi_{\mathrm{V}}^T$ are called the left and right eigenvectors, respectively. Also, $\Sigma$ is a diagonal matrix containing the eigenvalues organized in a descending order. Therefore, each element of the $\Sigma$ matrix demonstrates the importance of each column of $\Phi_{\mathrm{U}}$. It should be noted that, the left eigenvector corresponds to the spatial coherent structures which are in other words, the modes of the flow. On the other hand, the right eigenvector indicates the evolution of each mode in time domain.

\subsection{Validation Study}

\begin{figure}
	\centering
	\includegraphics[scale=.55]{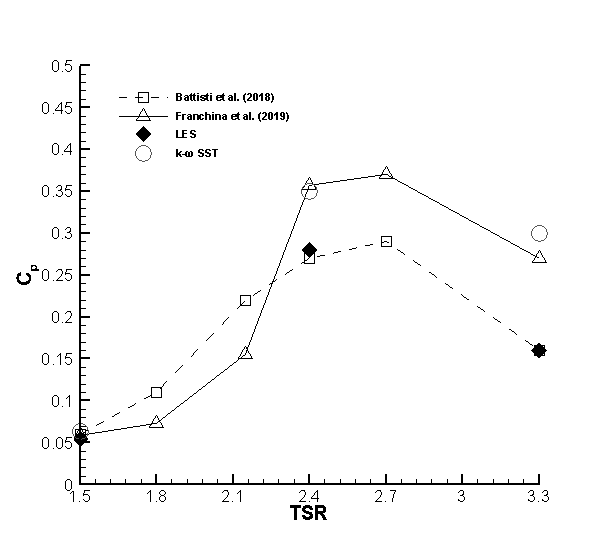}
	\caption{Power coefficient validation with respect to the experimental work of \citep{Battisti_2018} and numerical results of \citep{Franchina_2019}}
	\label{FIG4}
\end{figure}

The RANS and LES simulations were performed at the peak power TSR 2.4 and the two extremes of off-design points i.e. TSRs 3.3 and 1.5. Following this, the results were extracted and the moment coefficient versus TSR were plotted in Fig. \ref{FIG4}. The results are in good agreement with that of the experimental study of Battisti et al. \citep{Battisti_2018} and the RANS numerical results of Franchina et al. \citep{Franchina_2019,Franchina_2020}. It should be noted that the difference between the two-dimensional RANS and the experimental results stems from the fact that in 2-D simulations the effects of tip vortices and struts are absent and the flow experiences lower drag and therefore drop in power consequently which is in agreement with the results of Franchina et al. \citep{Franchina_2019,Franchina_2020}. Moreover, the main discrepancy between LES and RANS occurs as a result of the dynamic stall characteristics (i.e. formation, growth, bursting/
shedding of the separation bubble) and trailing-edge vortex as well as their interaction specially at low tip speed ratios. More details of RANS and SRS methods' applicability for simulating VAWTs can be found in \citep{Franchina_2019,Franchina_2020}. 

\begin{figure}
	\centering
	\includegraphics[scale=.5]{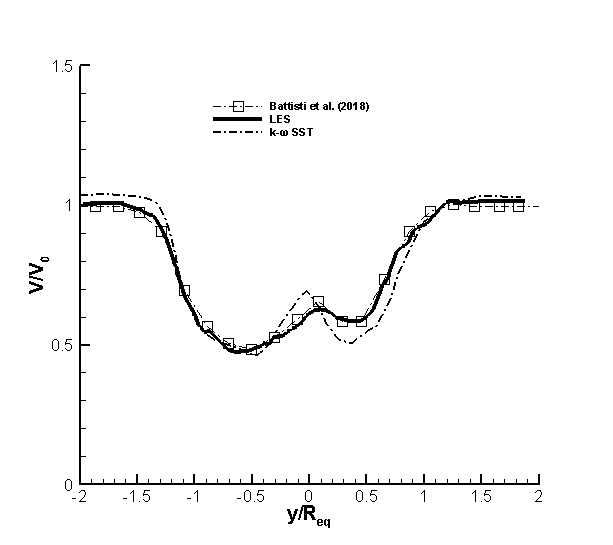}
	\includegraphics[scale=.5]{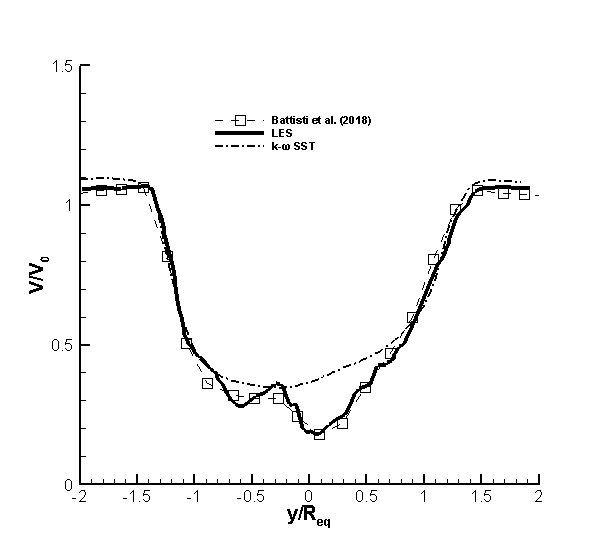}
	\includegraphics[scale=.5]{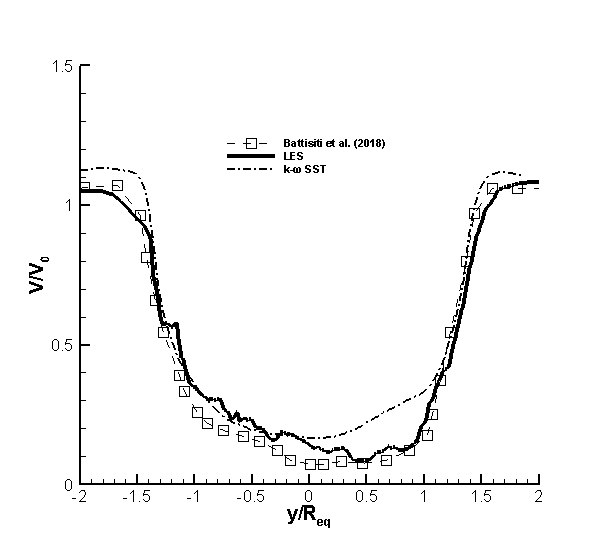}
	\caption{Validation in the wake region located at 3.5D behind the rotor with respect to the experimental results of \citep{Battisti_2018} a) TSR1.5, b) TSR2.4 and c) TSR3.3 }
	\label{FIG5}
\end{figure}

As the main purpose of the current study is to investigate the wake characteristics behind the VAWT, the velocity profile at the wake region located 3.5D behind the rotor was calculated and the results were compared against that of the experimental work of \citep{Battisti_2018} in Fig. \ref{FIG5}. Both RANS and LES approaches exhibited good agreement with the experimental data for the time-averaged velocity. However, the main shortcoming of RANS model, i.e. inability to capture unsteady and inherent fluctuations of the flow, resulted in a smoother velocity profile compared to that of LES \citep{lin2021cfd}.

\section{Results and Discussion}

\begin{figure}
	\centering
	\includegraphics[scale=.15]{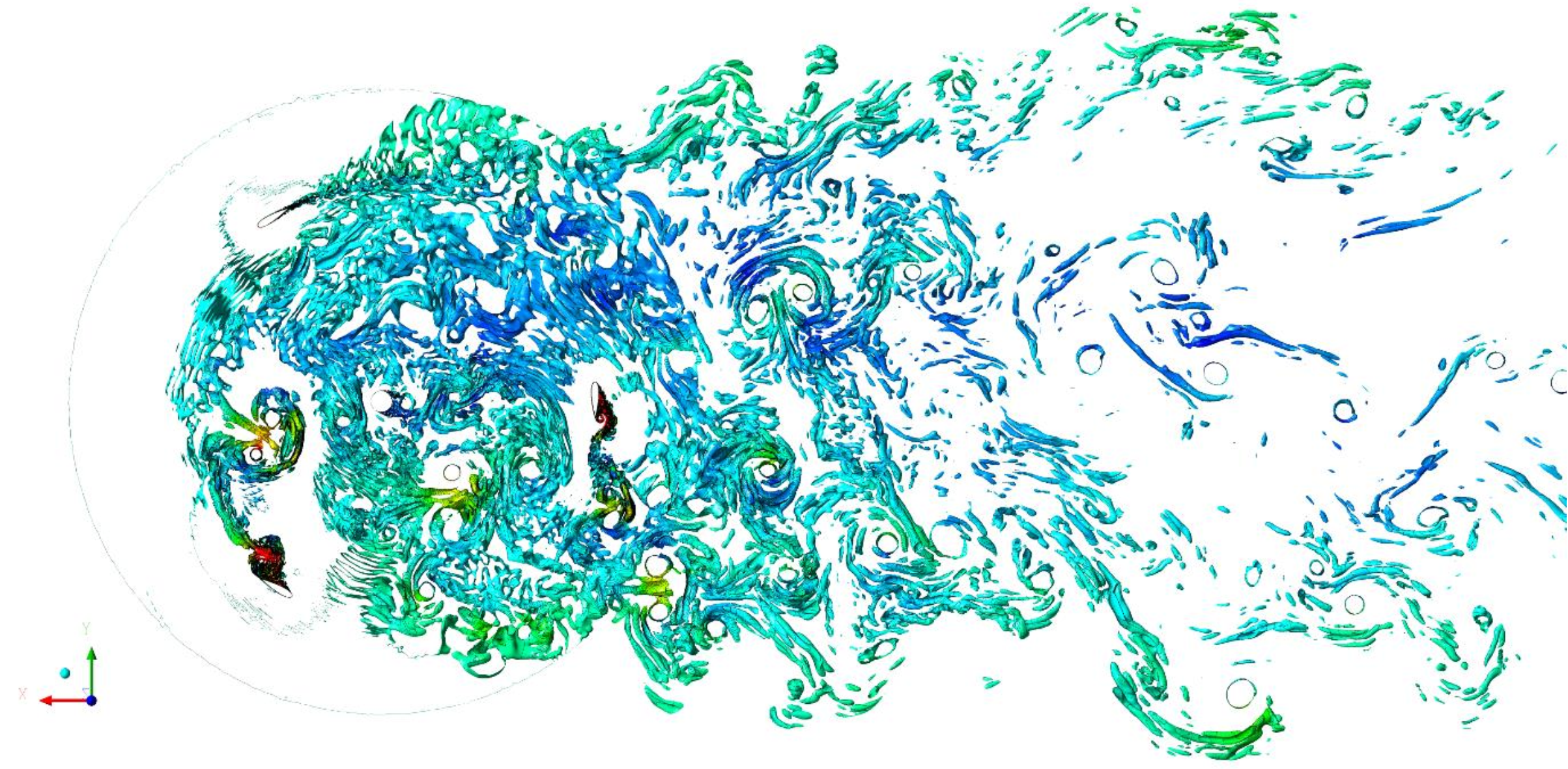}
	\includegraphics[scale=.06]{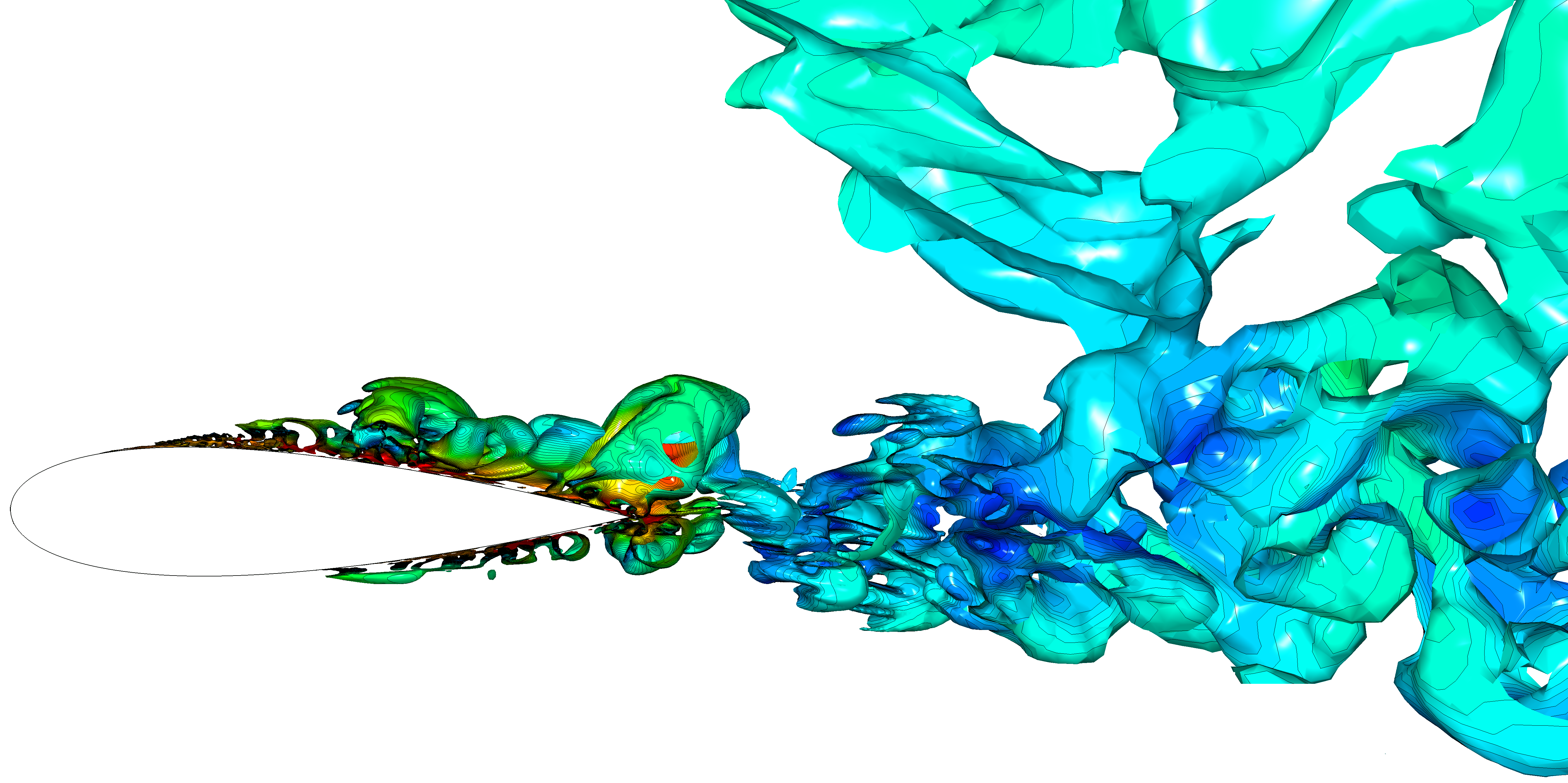}
	\includegraphics[scale=.06]{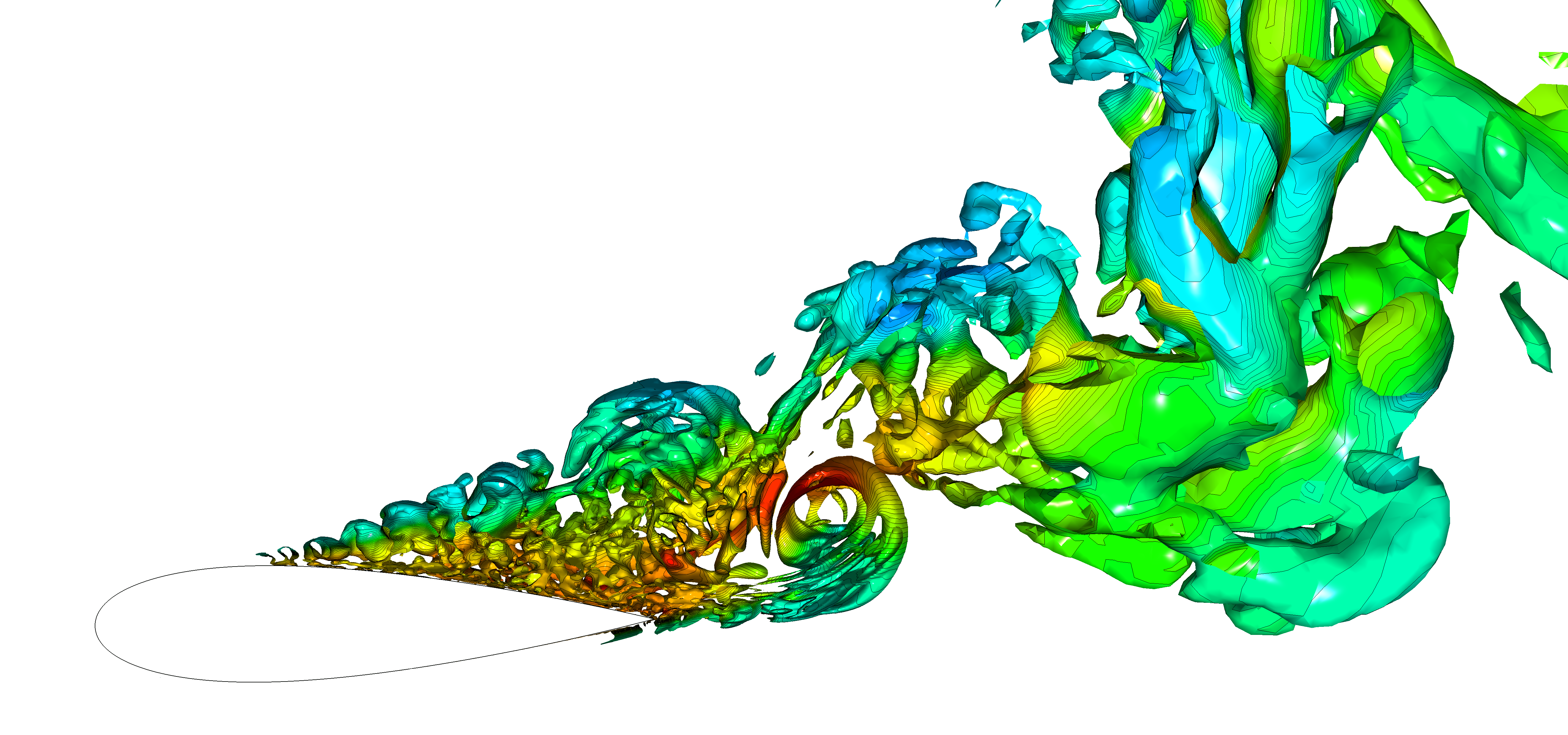}
	\includegraphics[scale=.06]{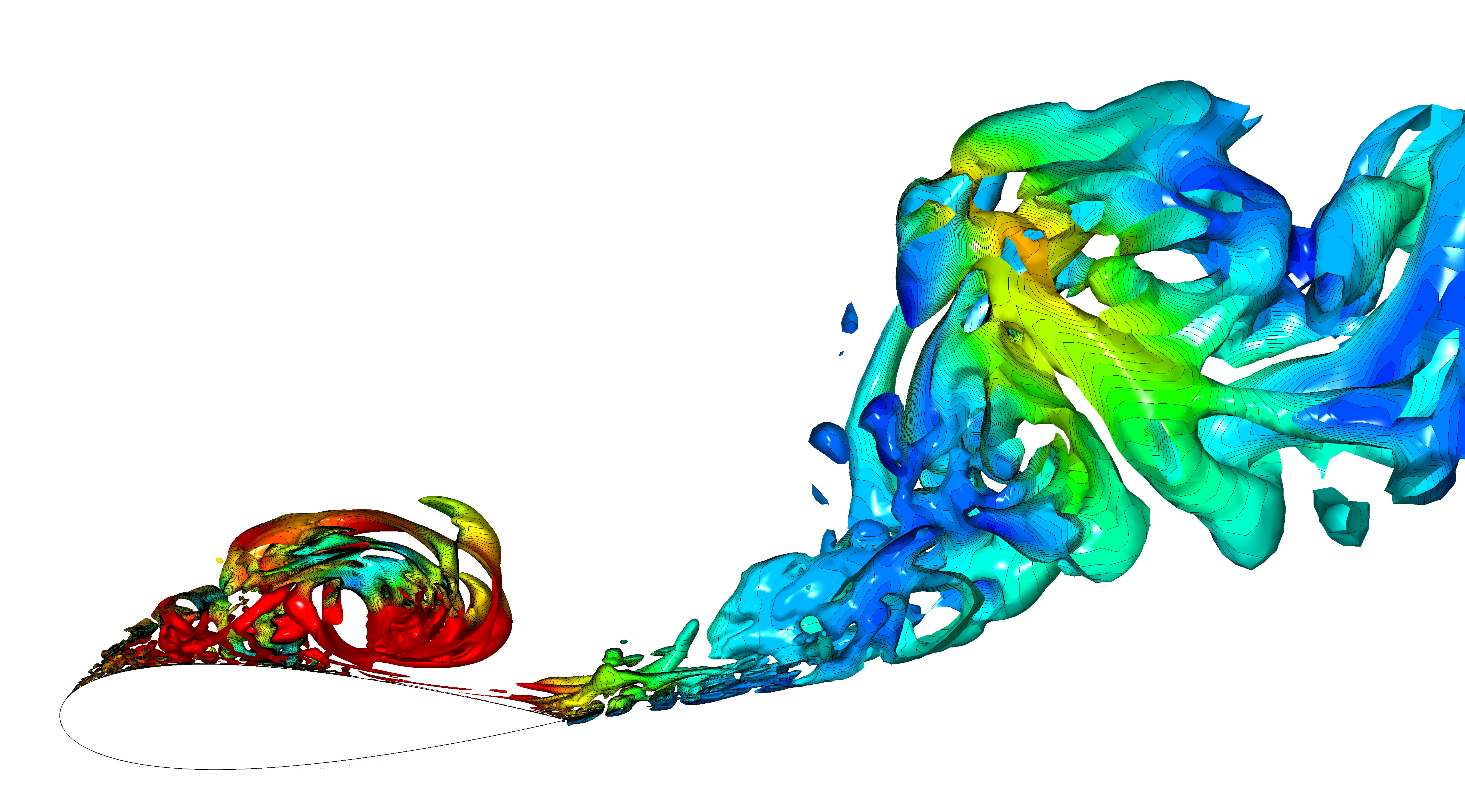}
	\includegraphics[scale=.1]{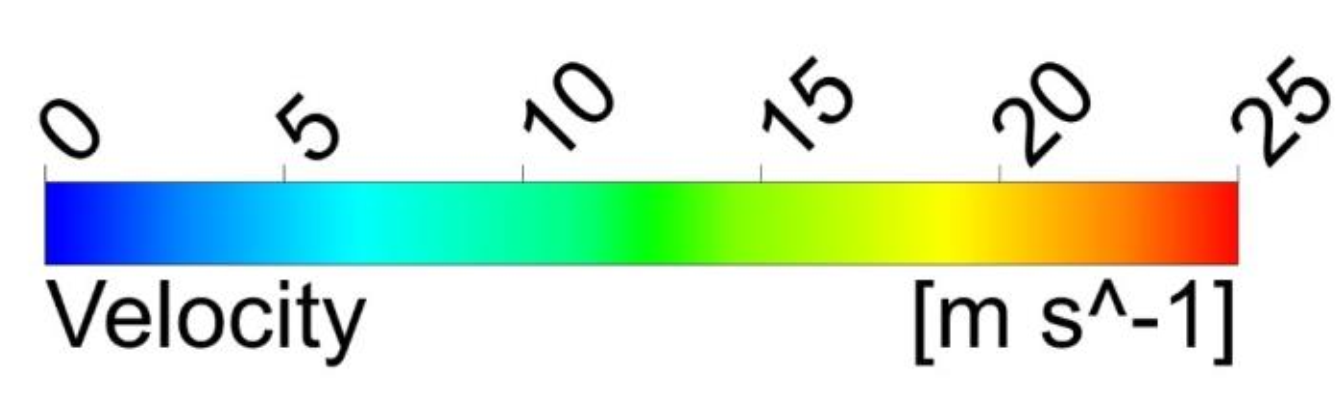}
	\caption{Vorticity contours of LES simulation of TSR2.4 at different rotation sections a) vorticity contour in the entire flow domain b) vortex Shedding at windward section c) vortex Shedding at upwind section d) vortex Shedding at leeward section}
	\label{FIG4_1}
\end{figure}

The numerical simulation was performed at different TSRs by both RANS and LES methods and the POD modes were extracted. In the first place, the flow physics and the structures formed due to the vortex shedding past the rotor will be briefly discussed. Fig. \ref{FIG4_1} depicts the vorticity contours colored by the velocity magnitude at different sections of the rotation domain namely windward, upwind, leeward and downwind whose definition and a thorough explanation on them can be found in Tescione et al. \citep{Tescione_2014}. As long as the blade sweeps the windward section (Fig. \ref{FIG4_1}-b), it sees a relatively low angle of attack which could be clearly observed for the one on the top of the rotating domain where the vortex shedding is lower compared to the other two airfoils. This is due to the fact that at this section the rotational velocity forms a considerably low angle with the absolute velocity vector resulting in a low angle of attack and blade loading and therefore less separation. However, as the blade progresses towards the upwind section Fig. \ref{FIG4_1}-c) the vortex shedding sees a rising trend due to the increased blade loading. As the blade approaches the leeward section (Fig. \ref{FIG4_1}-d), the separation reaches its height and having entered the downwind part, the separation and the consequent vortex shedding phenomenon abates.    
 
\begin{figure}
	\centering
	\includegraphics[scale=.5]{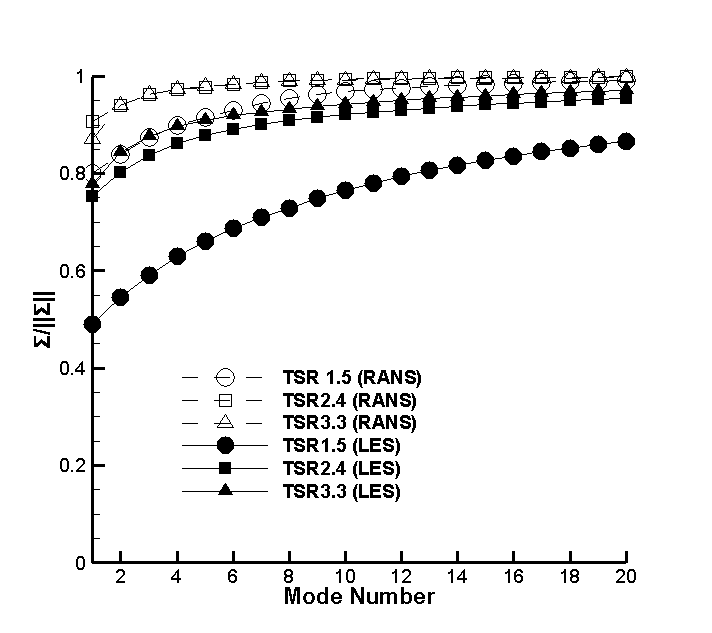}
	\includegraphics[scale=.5]{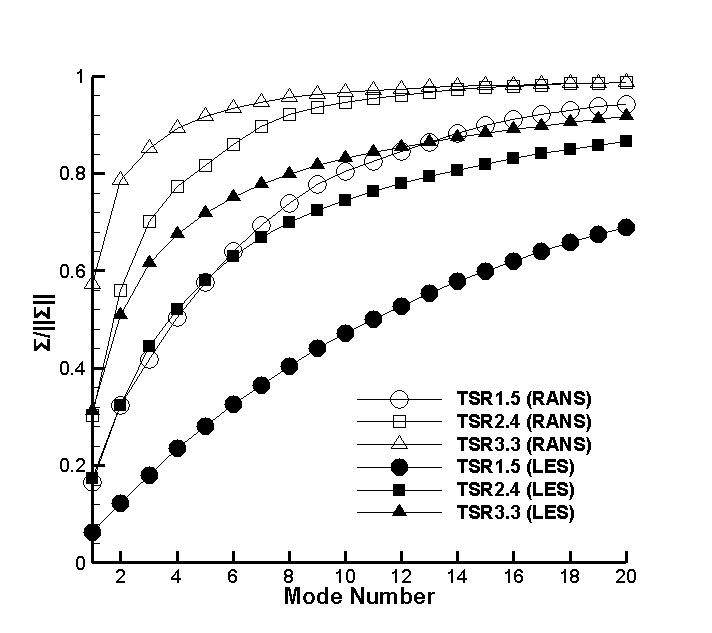}
	\caption{Cumulative POD energy curves a) streamwise component of velocity b) transversal component of velocity}
	\label{FIG4_2}
\end{figure}

The cumulative modal energy of the RANS and LES results have been presented in Fig. \ref{FIG4_2}. Regarding the RANS results, the cumulative modal energy of the streamwise component of velocity proves that except for the first mode, the modal energy of the TSR2.4 and TSR 3.3 fall on top of each other. However, TSR1.5 reaches the maximum amount of the accumulative modal energy at significantly higher mode numbers suggesting that more flow structures could be excepted at lower TSRs due to the higher blade loading and high relative angle of attack which is an intrinsic characteristic of lower TSRs. Also, as shown in Fig. \ref{FIG4_2}-b, the number of effective modes increases considerably with the decrease in TSR. As a result, it is found that the vortex shedding due to separation at high relative angles of attack leads to more effective modes for transversal component of velocity rather than the streamwise component of velocity. Moreover, comparing to the streamwise component of velocity, it is observed that the largest amount of RANS modal energy accumulation for the transversal component of velocity belongs to the highest TSR, while the peak power (TSR 2.4) has the highest amount for the streamwise velocity component with respect to the other cases. Therefore, the transversal modes mainly correspond to flow structures stemming from flow separation and entrainment which is a well-known characteristic of lower TSRs. Therefore, it is observed that the number of convective modes (effective modes) is proportional to the wind velocity for the transversal component modes. In addition, due to the striking similar behavior of the cumulative energy between TSR3.3 and TSR2.4, analogous spatial modal behavior could be expected. It should be noted that this observation could be better understood by plotting the spatial POD modes.

However, as for the LES results, the modes of TSR2.4 and TSR3.3 show a more noticeable difference with respect to the RANS modes, as shown in Fig. \ref{FIG4_2}-a. This is due to the fact that the LES resolves more structures compared to RANS simulation leading to more convective modes. Therefore, the observation previously made on RANS modes is confirmed suggesting that the lower TSRs are correspondent to more wake shedding which is in agreement with the aerodynamics of VAWTs. Also, it should be noted that unlike the RANS cumulative modes, where the rated TSR demonstrated the maximum energy at the first mode, the highest energy at the same mode number is observed at the highest TSR for the LES results. This observation suggests that the RANS fails to properly compute the energy of the modes at higher TSRs. In other words, the effective structures at higher TSRs cannot be distinctively captured in RANS POD modes. Moving on to Fig. \ref{FIG4_2}-b, the cumulative modal energy behavior of the transversal component of velocity reveal that both RANS and LES still follow the same trend except that the LES modes, as previously stated, demonstrate considerably a lower amount of modal energy at each mode number owing to the capability of LES method to resolve considerably more flow structures.

\begin{figure}
	\centering
	\includegraphics[scale=.7]{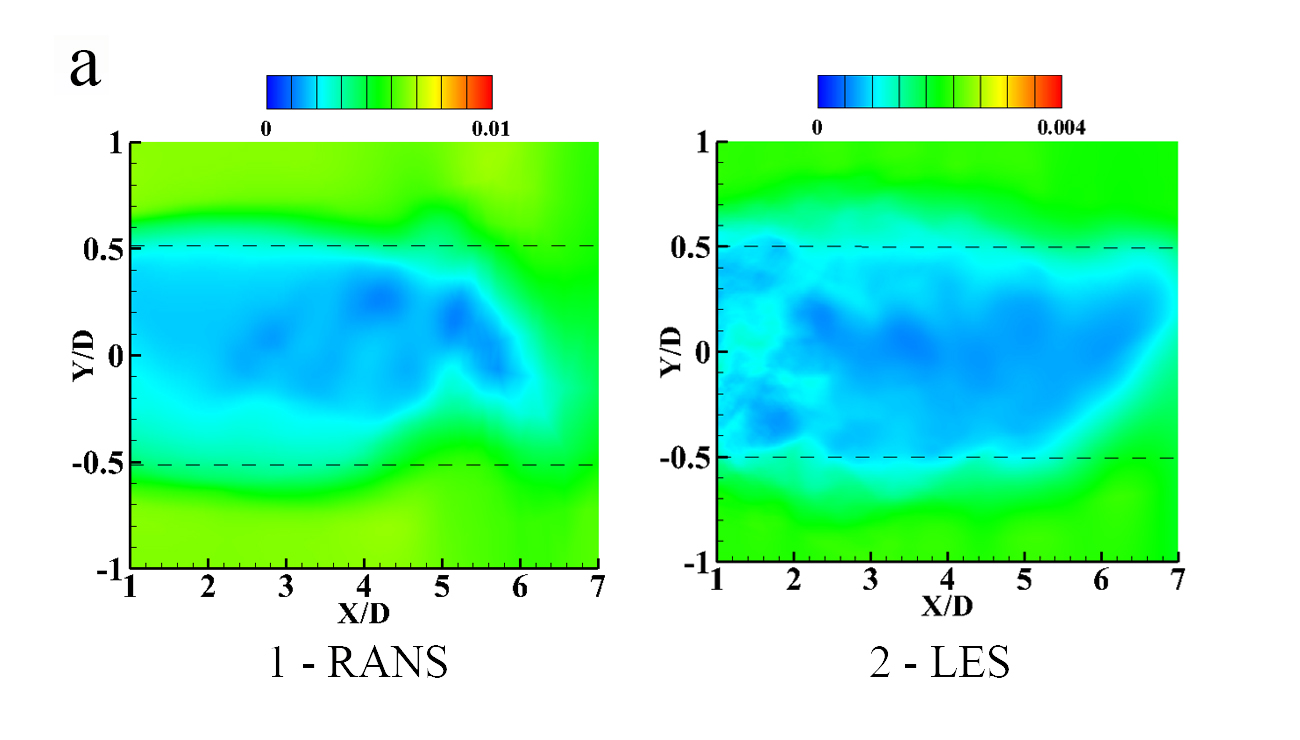}
	\includegraphics[scale=.7]{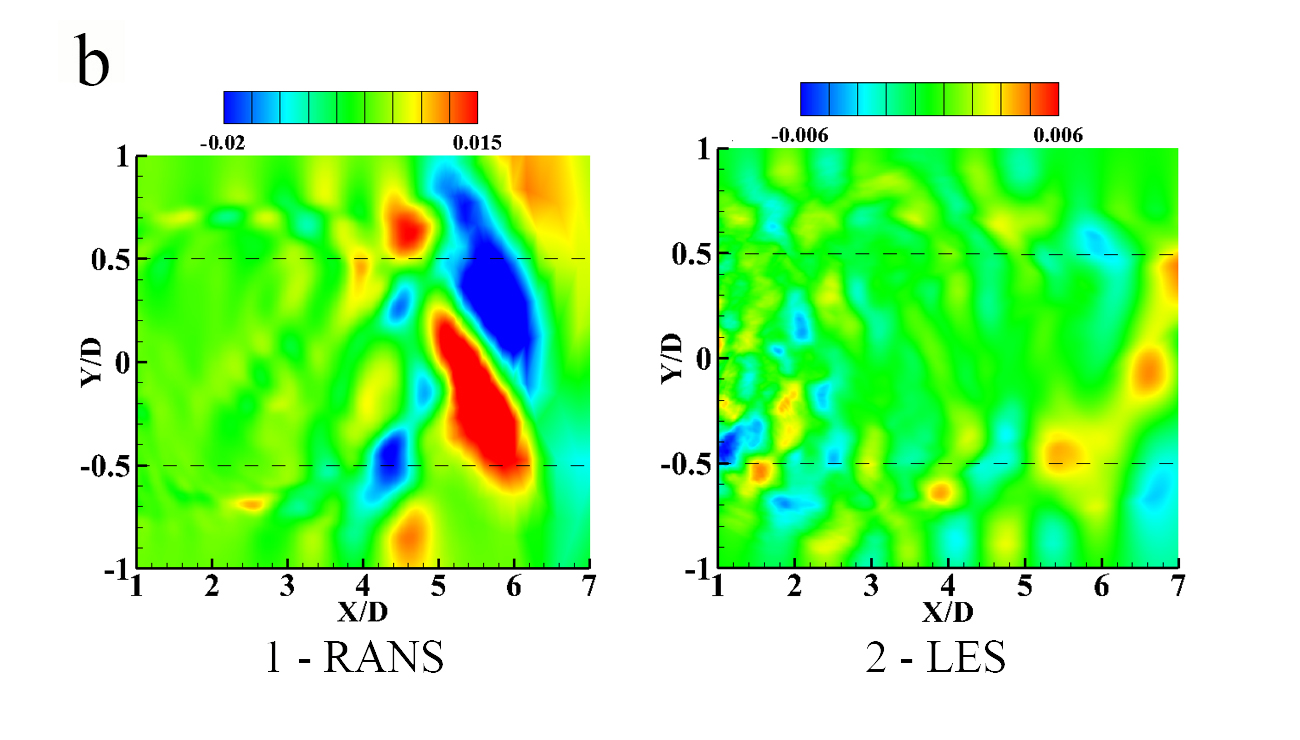}
	\includegraphics[scale=.7]{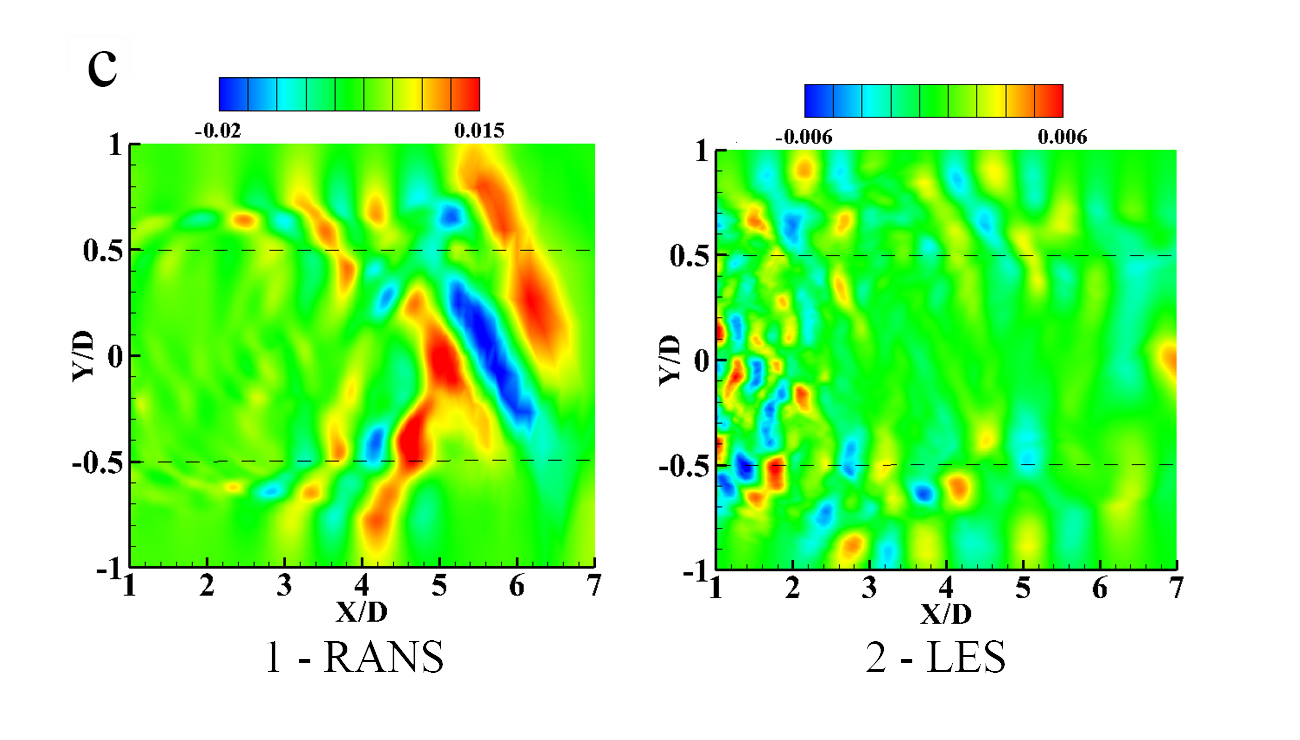}
	\includegraphics[scale=.7]{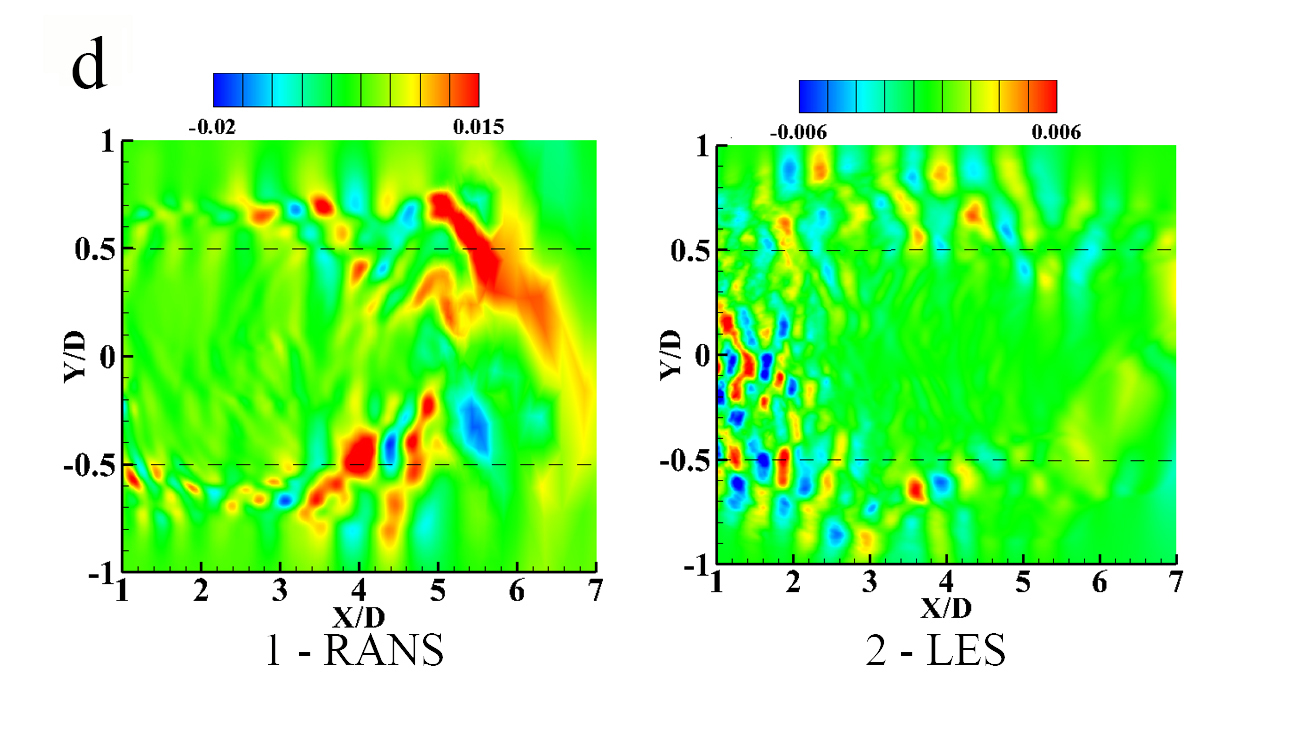}
	\includegraphics[scale=.7]{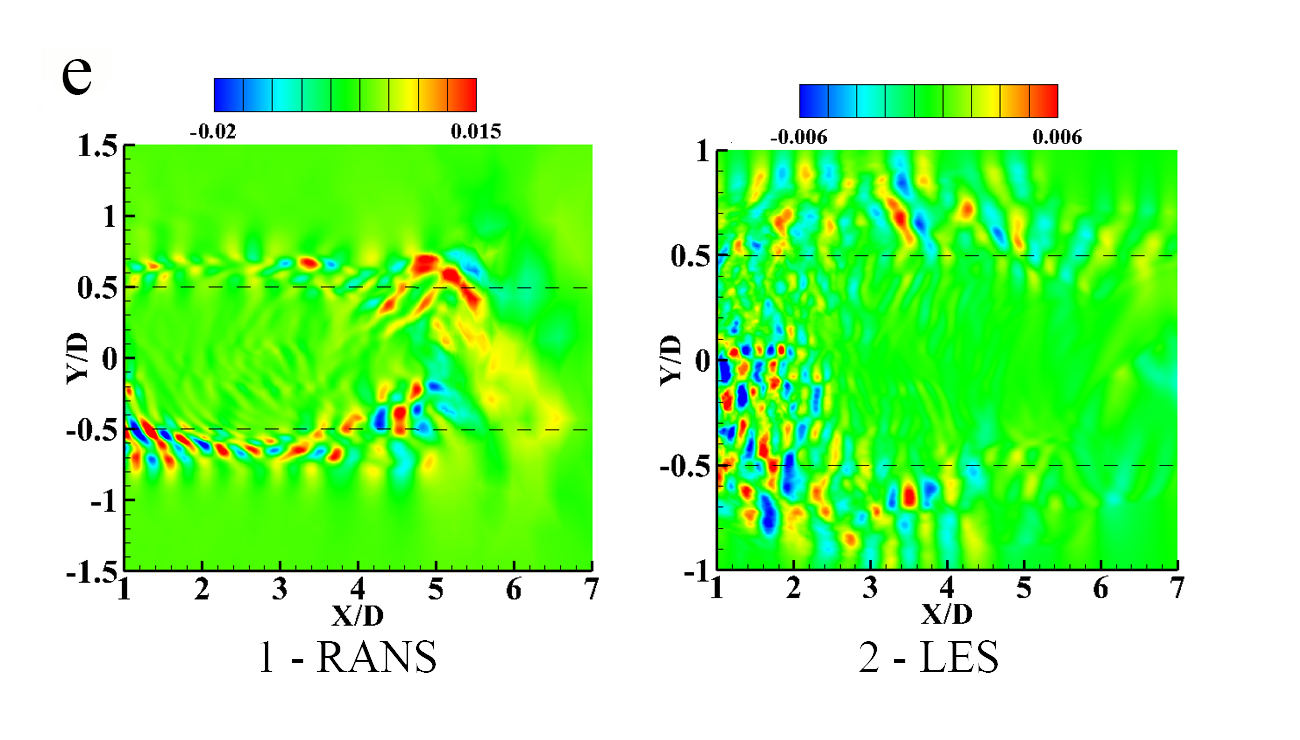}
	\caption{POD modes for streamwise component of velocity at TSR2.4 (left RANS and right LES)}
	\label{FIG4_3}
\end{figure}

\begin{figure}
	\centering
	\includegraphics[scale=.7]{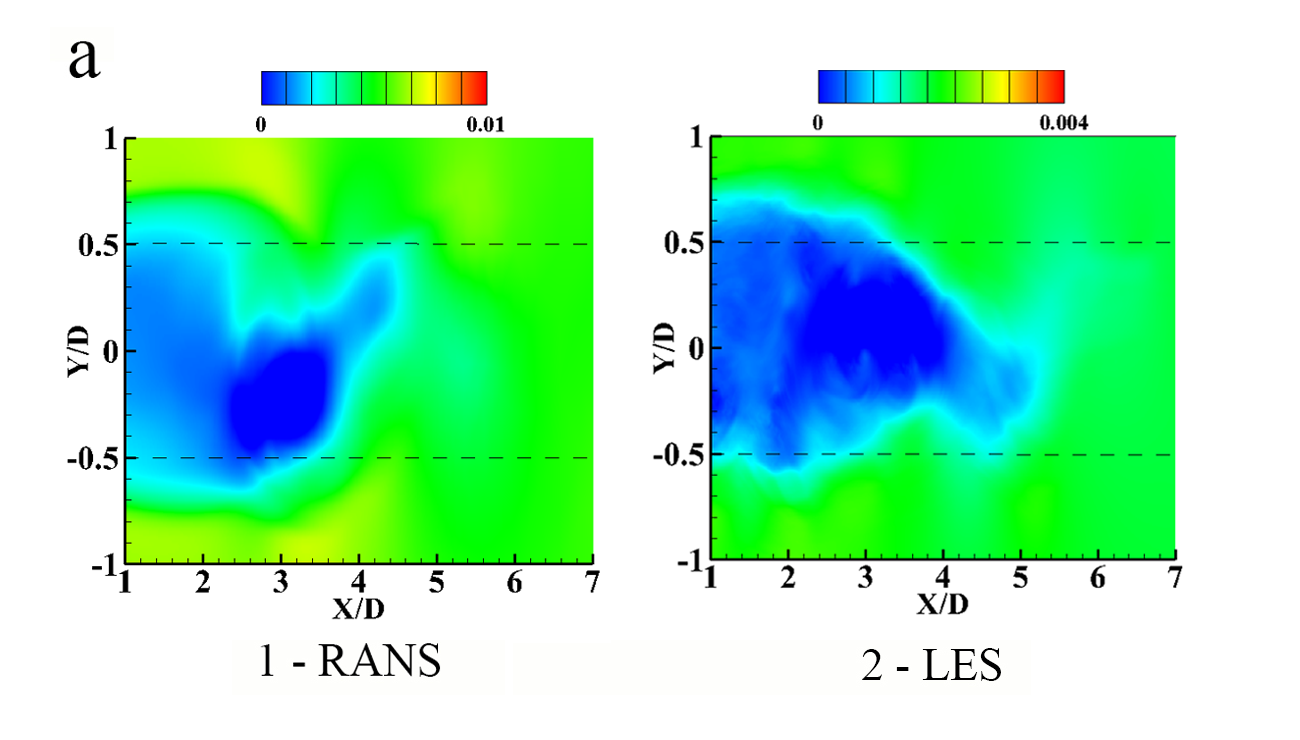}
	\includegraphics[scale=.7]{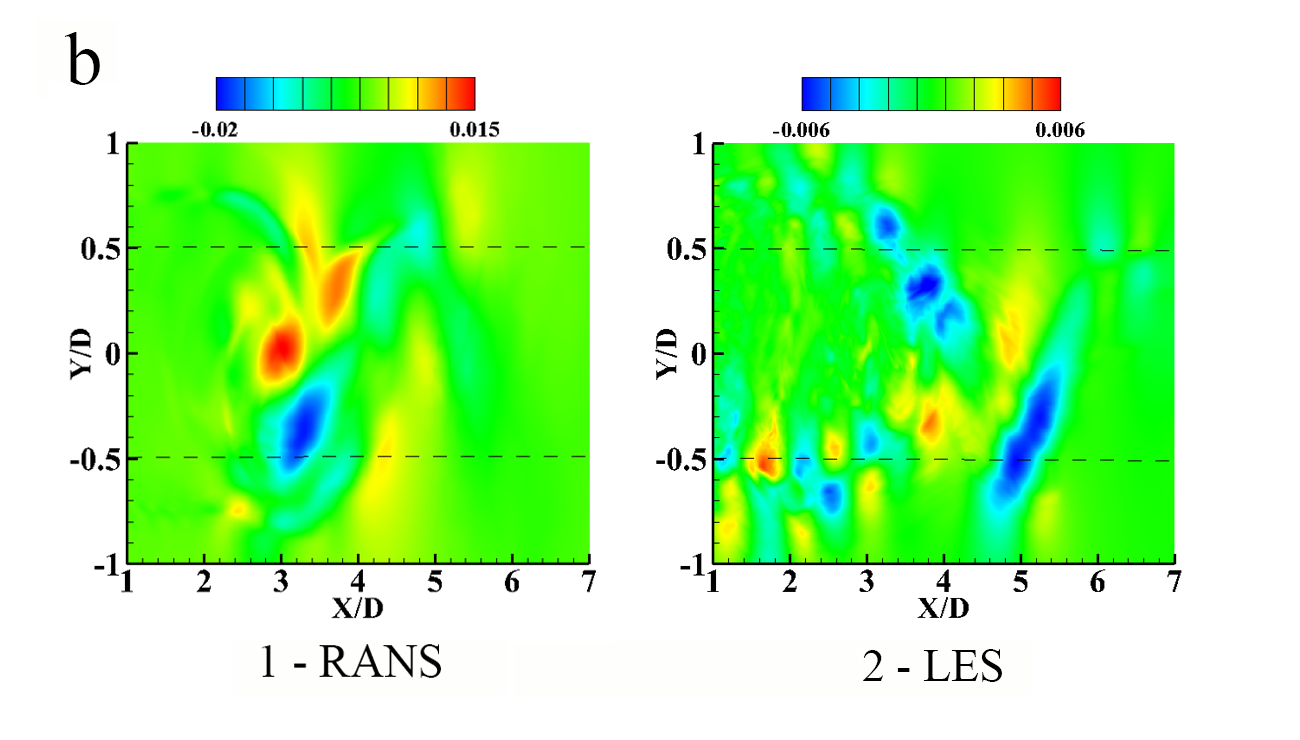}
	\includegraphics[scale=.7]{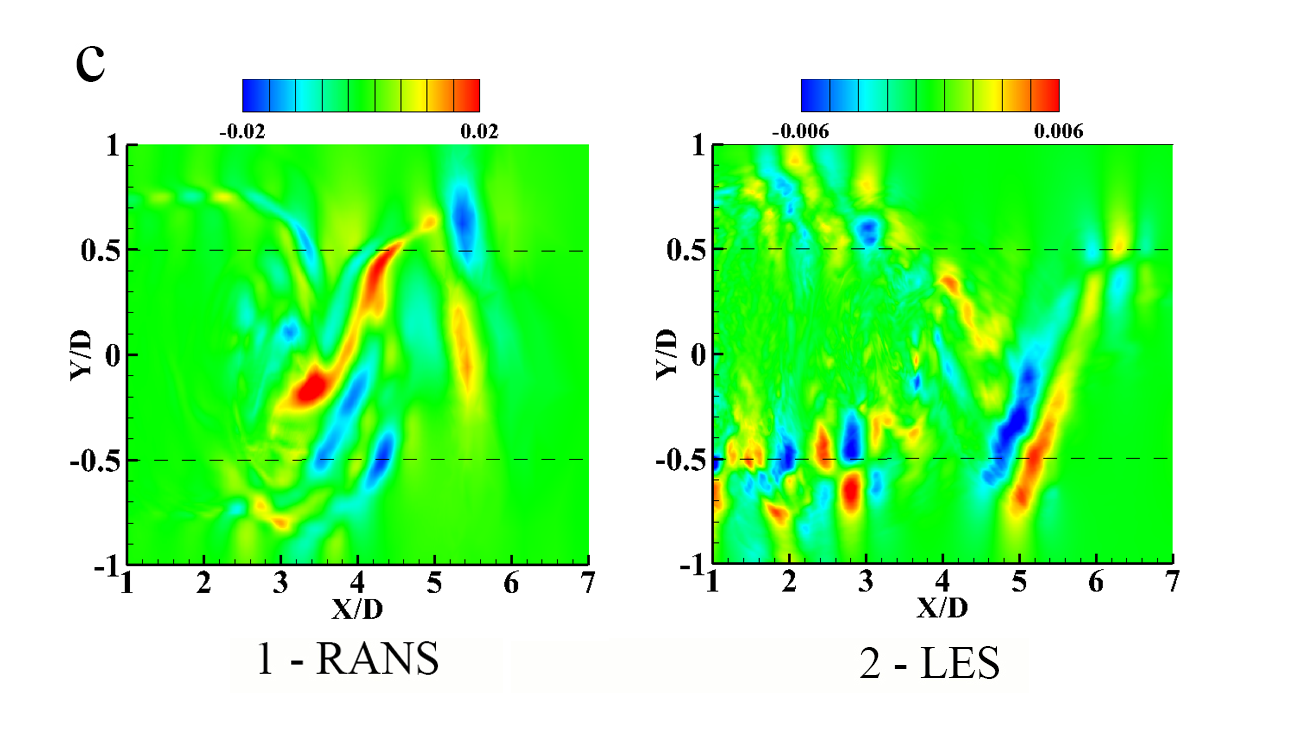}
	\includegraphics[scale=.7]{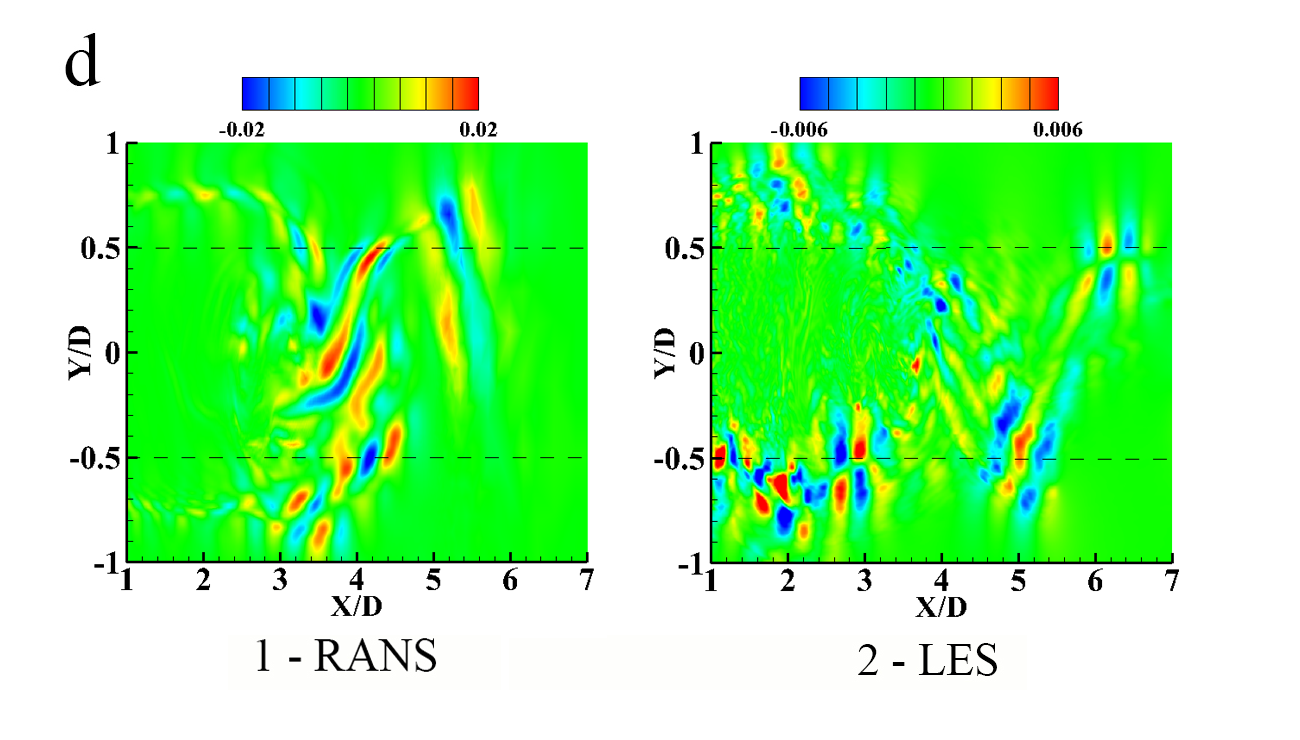}
	\includegraphics[scale=.7]{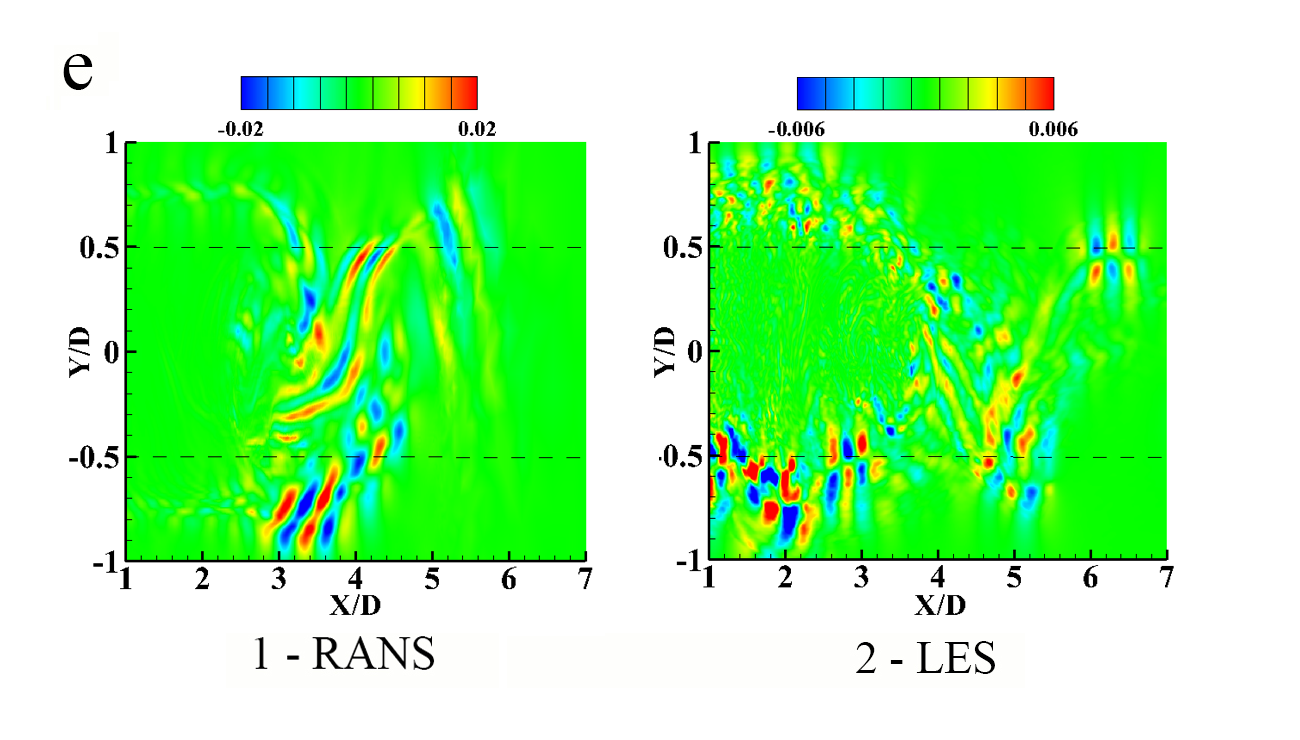}
	\caption{POD modes for streamwise component of velocity at TSR3.3 (left RANS and right LES)}
	\label{FIG4_4}
\end{figure}

\begin{figure}
	\centering
	\includegraphics[scale=.7]{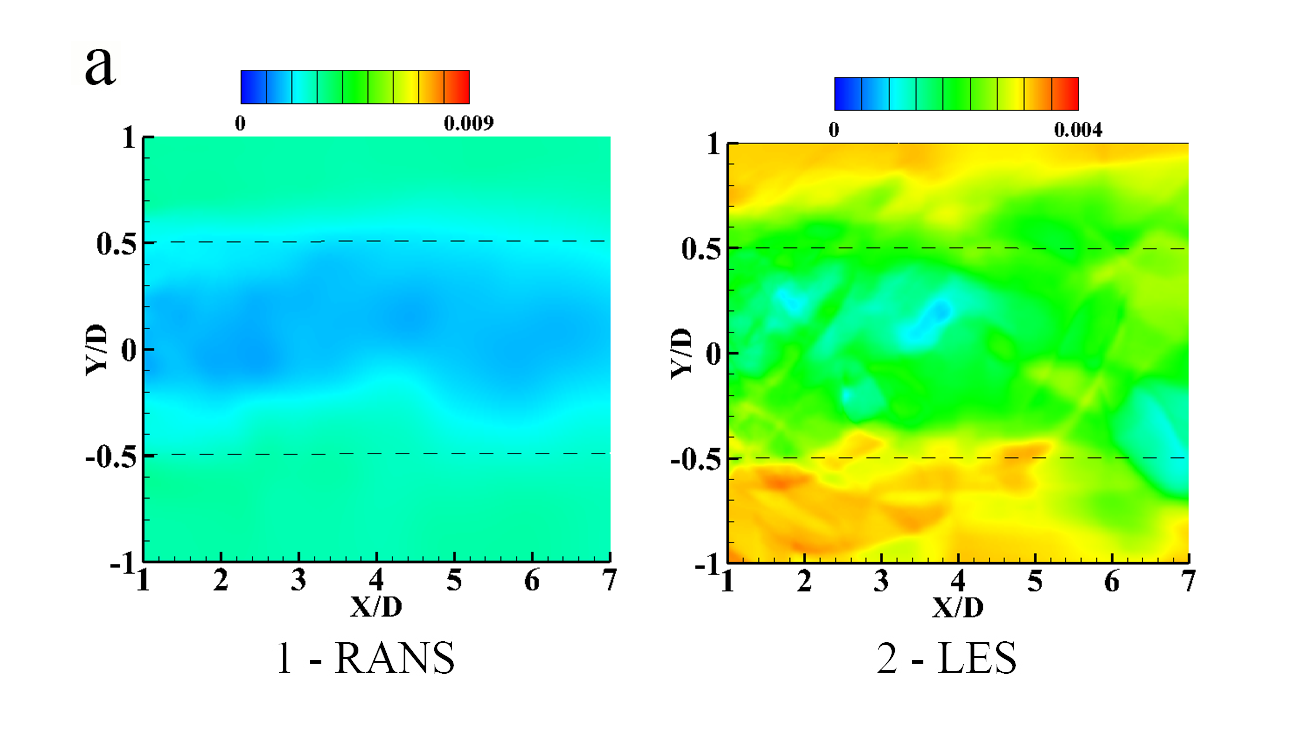}
	\includegraphics[scale=.7]{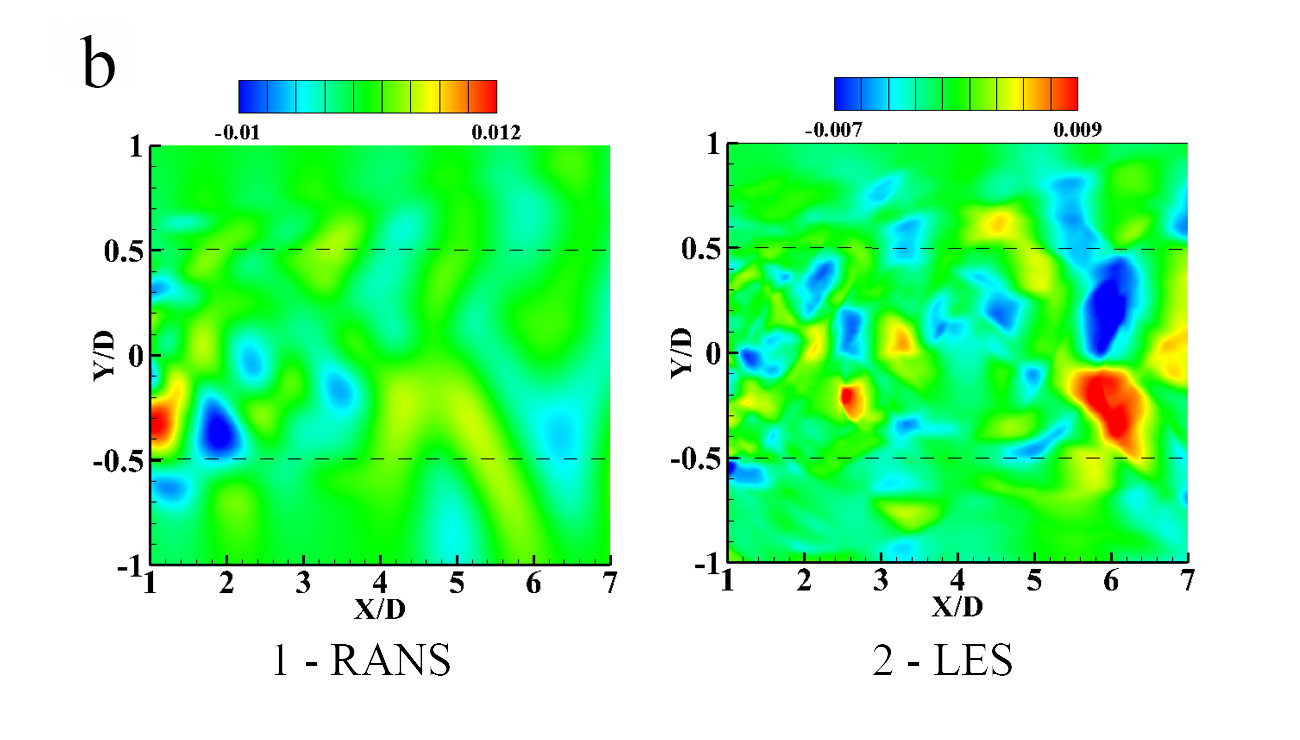}
	\includegraphics[scale=.7]{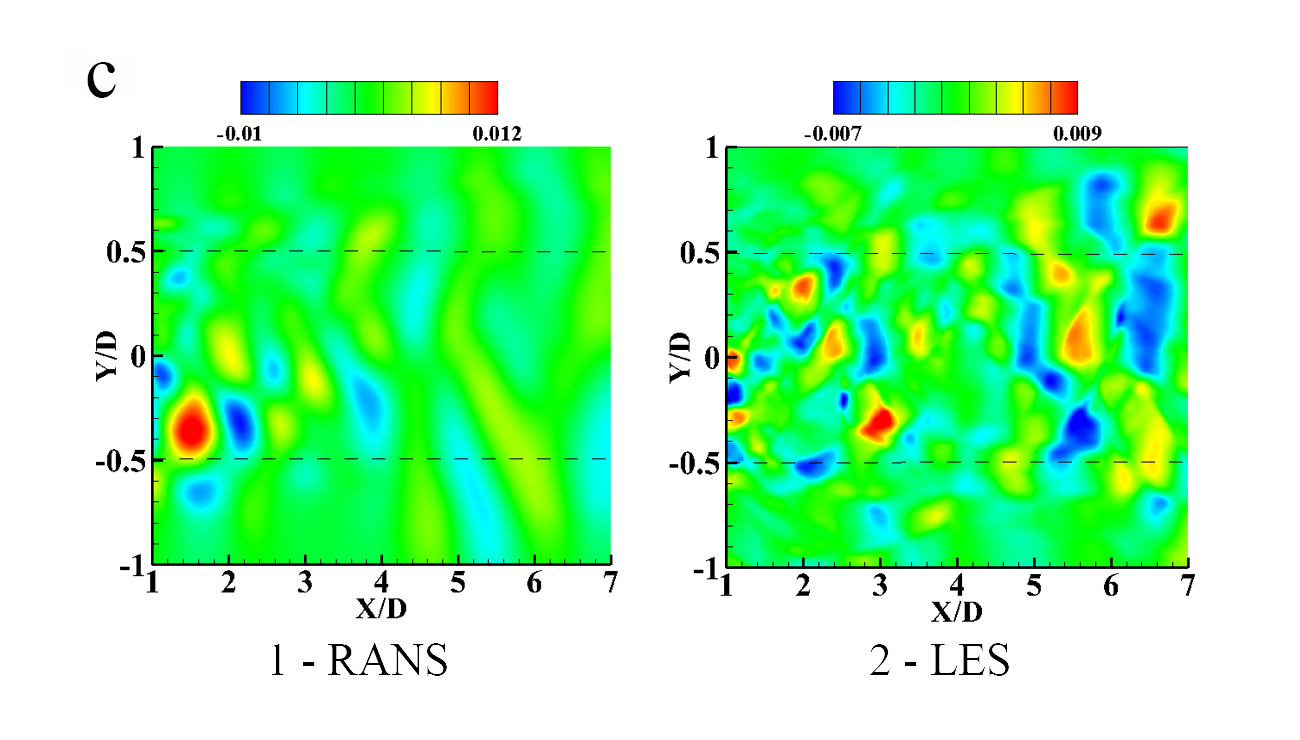}
	\includegraphics[scale=.7]{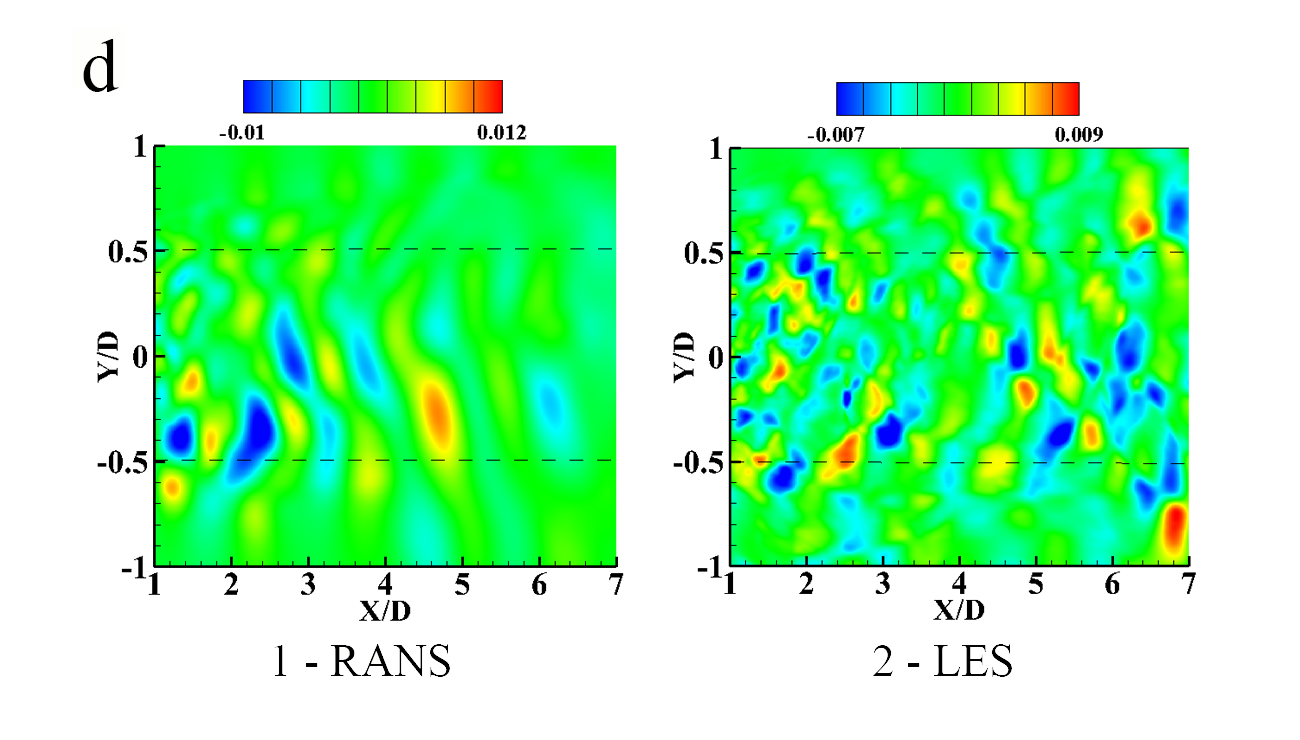}
	\includegraphics[scale=.7]{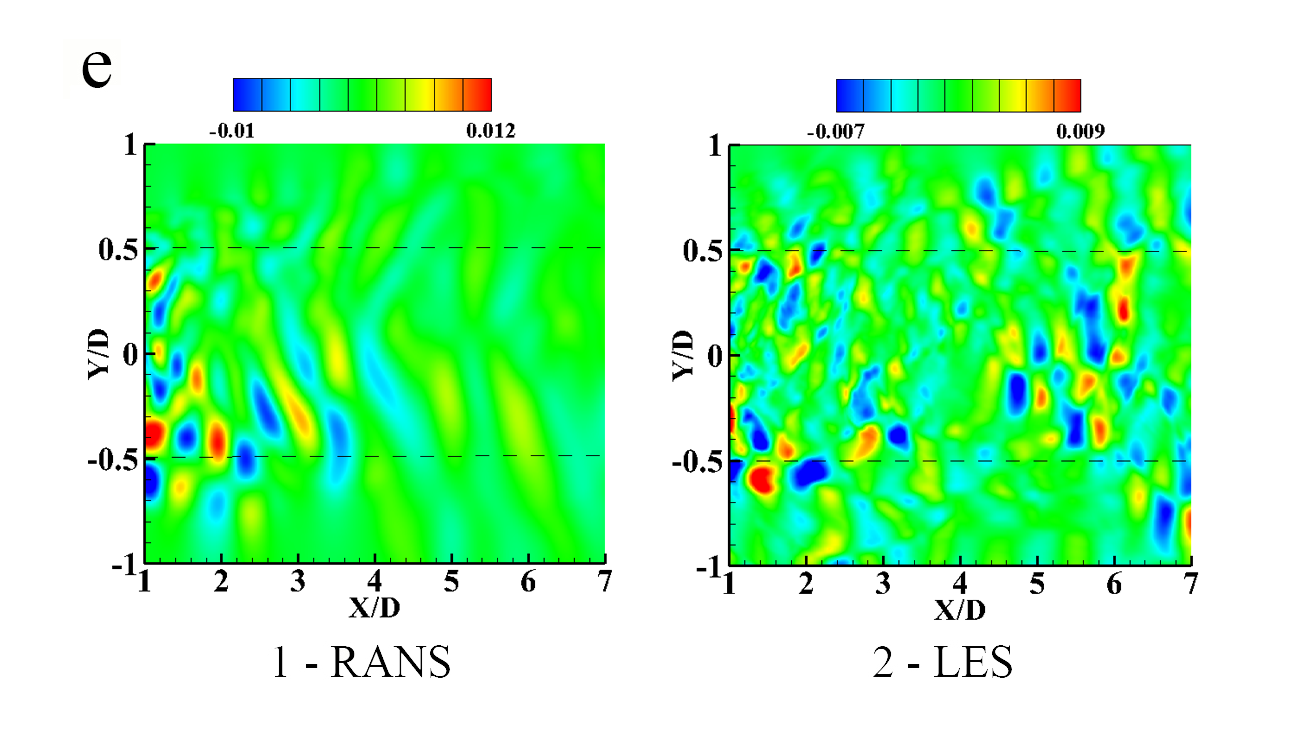}
	\caption{POD modes for streamwise component of velocity at TSR1.5 (left RANS and right LES)}
	\label{FIG4_5}
\end{figure}

Presented in Fig. \ref{FIG4_3},  \ref{FIG4_4} and \ref{FIG4_5} are the spatial POD modes at TSR2.4, TSR3.3 and TSR1.5 respectively where the distance between the dashed lines indicates the diameter of the wind turbine as they have been placed at the two extremes of the rotor. Also, the windward and leeward sections of the wind turbine have been located at the upper and lower boundaries of the POD box, respectively. As for the mode analysis of RANS simulations, the first mode shown in Fig. \ref{FIG4_3}-a-1, demonstrates the mean flow behind the rotor which could be interpreted as the efficiency of the turbine in energy extraction when compared to same mode number at TSRs 3.3 and 1.5 presented respectively in Fig. \ref{FIG4_4}-a-1 and Fig. \ref{FIG4_5}-a-1. It is evident in Fig. \ref{FIG4_3}-a-1 that the area affected by the presence of the wind turbine rotor, i.e. wake region, begins to collapse after about 6D, while for TSR3.3 and 1.5 the wake region collapse is observed at a distance considerably closer and far from to the rotor, respectively (Fig. \ref{FIG4_4}-a-1 and Fig. \ref{FIG4_5}-a-1). It should be noted that this is in agreement with our knowledge of physics of fluid flow past a VAWT. To illustrate, as the higher TSRs correspond to lower blade loading and less flow momentum, the fluid flow past the VAWT rotor experiences a wake region collapse at a closer distance to the rotor compared to the lower TSRs. On the other hand, as the lower TSRs are associated with higher blade loadings and fluid momentum, a larger area behind the rotor is affected by the energy extraction from the turbine.

Higher mode numbers reveal more vortical structures in smaller dimensions in both wake region and far field area. For instance, the second spatial mode shown in Fig. \ref{FIG4_3}-b-1, demonstrates the main flow structures at the point of wake region collapse. It should be noted that the same flow structures are observed at the same location in Fig. \ref{FIG4_4}-b-1. In addition, with the further rise in the mode number, the flow structures, in the wake region collapse area, decay (Fig. \ref{FIG4_3}-c-1) and the rolled-up flow structures at the edges of the wake region (shown in dashed lines) become more apparent as seen in Fig. \ref{FIG4_3}-d-1 and Fig. \ref{FIG4_3}-e-1 for the fourth and the fifth modes respectively. The shear flow at the edges of the wake region is due to the velocity difference between the wake region and the far field, as it is evident in the first mode, which results in instabilities leading to the formation of symmetrical rolled-up structures. Further, it is evident in Fig. \ref{FIG4_3}-d-1 and Fig. \ref{FIG4_3}-e-1 that the structures at the border of the leeward section is more pronounced than that of the windward section. This is attributed to the fact that the vortex shedding phenomenon is more considerable at leeward section due to the higher relative angle of attack the blade sees with respect to the windward section. Also, the same modal behavior is observed at TSR3.3 as it is evident in Fig. \ref{FIG4_4}-d-1 and Fig. \ref{FIG4_4}-e-1. Therefore, the modal behaviors at TSR2.4 and TSR3.3 are quite similar which was also shown in the modal cumulative energy plots. However, this similarity appears to be less noticeable at TSR1.5. Unlike TSR2.4 and TSR3.3 in which the wake region spans an area with almost the diameter of the wind turbine rotor, the wake region in TSR1.5 encompasses mostly the windward and parts of the downwind sections. This observation is completely in agreement with the fact that lower TSRs are associated with significantly lower amount of power coefficient as a less extensive span of the downstream flow is affected by the rotor presence compared to TSR2.4 and TSR3.3. Also, the spatial modes at TSR1.5 are prone to more vortex shedding due to the higher blade loading and angle of attack. Moreover, flow structures emanating from flow separation over the wind turbine blades are observed in all spatial modes of streamwise POD modes at TSR1.5. It should be noted that due to the fact that the flow undergoes lower velocity drop due to work extraction at TSR1.5 and consequently lower velocity difference between the wake region and the far field area, the symmetrical rolled-up structures are apparent to a lower extent at the edges of the wake region.

Generally, the second and the subsequent modes correspond to the decay of the wake region collapse structures and the structures created as the result of the velocity gradient between the wake region and the far field area. Thus, the second mode depicted in Fig. \ref{FIG4_3}-b-1 shows the structures emanating from the wake region collapse which is the same characteristics of all the TSRs at this mode number for the streamwise component of velocity. Moreover, with the increase in the mode number, the flow structures at the wake region collapse point begin to disappear and instead the vortical structures stemming from the velocity gradient at the edge of the wake region become more pronounced which are evident in Fig. \ref{FIG4_3}-d-1 and Fig. \ref{FIG4_3}-e-1. 

In Fig. \ref{FIG4_3}-a-2, Fig. \ref{FIG4_4}-a-2 and Fig. \ref{FIG4_5}-a-2 the first five dominant modes of LES results for TSRs2.4, 3.3 and 1.5 have been respectively presented. The first modes for all TSRs bear a resemblance with those of the RANS results which is due to the non-connectiveness nature of the first mode. In addition, comparing the first spatial modes of RANS and LES results, it reveals that both the numerical methods manage to predict the wake region collapse at approximately the same distance except that the LES suggests a further collapse position. However, the discrepancies between the RANS and LES modes begin to be more evident with the increase in the mode number which is in line with the expectations of the cumulative modes plots of both numerical methods. In this regard, as for the second LES mode at TSR2.4 presented in Fig. \ref{FIG4_3}-b-2, it is evident that in addition to the flow structures corresponding to wake region collapse which exists in the RANS mode as well, the wake shed emanating from the rotor is also observed which are bound to the boundaries of the POD box. Therefore, the LES results, contrary to the RANS modes, suggest more interpretable fluid phenomena in the wake region as expected due to the resolving the major part of the flow field. As for the third LES mode of the streamwise component of velocity at TSR2.4 presented in Fig. \ref{FIG4_3}-c-2 it is observed that firstly, unlike the third RANS mode, the structures at the edges of the POD box appear more pronouncedly near the rotor. In addition, the structures emanating from the leeward section are more noticeable which is due to the relative angle of attack at this part of rotation. The second feature of the third LES mode which is worth mentioning is that the structures shedding from the rotor in the preceding mode have decayed into smaller ones. The fourth LES mode of the streamwise component, Fig. \ref{FIG4_3}-d-2, suggests that in addition to the convective structures from the preceding modes which are resultants of the turbulence decay of larger structures, the symmetrical structures near the edges of the POD box appear more prominently near the windward section. In addition, comparing the fourth RANS and LES modes it is observed that while the RANS results demonstrate the main concentration of the structures to be near the wake region collapse area, the LES mode suggests the most pronounced area of entrainment to be at the beginning of the domain. In regards with the TSR3.3 modes presented in Fig. \ref{FIG4_3}-b, more similarity between the RANS and LES modes are observed. This is in agreement with the cumulative energy in Fig. \ref{FIG4_2}-a where the RANS and LES modes for TSR3.3 show closer values compared to the TSR2.4. As a result, the spatial modal behavior for RANS and LES at TSR3.3 follow a more analogous trend. To illustrate, the second LES mode of TSR3.3 in Fig. \ref{FIG4_4}-b-2 demonstrates the structures in the wake region collapse position which is the main spatial modal behavior of the second RANS mode inFig. \ref{FIG4_4}-b-1. As for the third LES mode at TSR3.3 (Fig. \ref{FIG4_4}-c-2), similar to the counterpart of the RANS mode as seen in Fig. \ref{FIG4_4}-c-1, the structures due to the mixing in the vicinity of wake region collapse and the edges of the POD box are observed, however, of different sizes. Furthermore, the same trend in spatial modal behavior is observed between LES and RANS in the subsequent modes presented in Fig. \ref{FIG4_4}-d and Fig. \ref{FIG4_4}-e.

Regarding the TSR1.5 LES modes, more discrepancies with respect to the TSRs2.4 and 3.3 are observed which is in agreement with the inference drawn from the RANS modes. In addition, it should be noted that more similarity between the RANS and the LES modes compared to the TSR2.4 can be found at TSR1.5 which was also discussed in cumulative curves for TSR3.3. Therefore, so far, we have come to this conclusion that the spatial modes at off-design points of operation can also be captured by RANS at a relatively acceptable accuracy. Also, it should be noted that unlike the spatial modes of the TSR3.3 the structures at TSR1.5 are mainly concentrated inside the wake region rather than the POD box edges which is due to the high blade loading and the lowest efficiency among the three TSRs. This trend is similarly followed in all the spatial LES and RANS modes of TSR1.5 as seen in Fig. \ref{FIG4_5}-c to Fig. \ref{FIG4_5}-e.

\begin{figure}
	\centering
	\includegraphics[scale=.7]{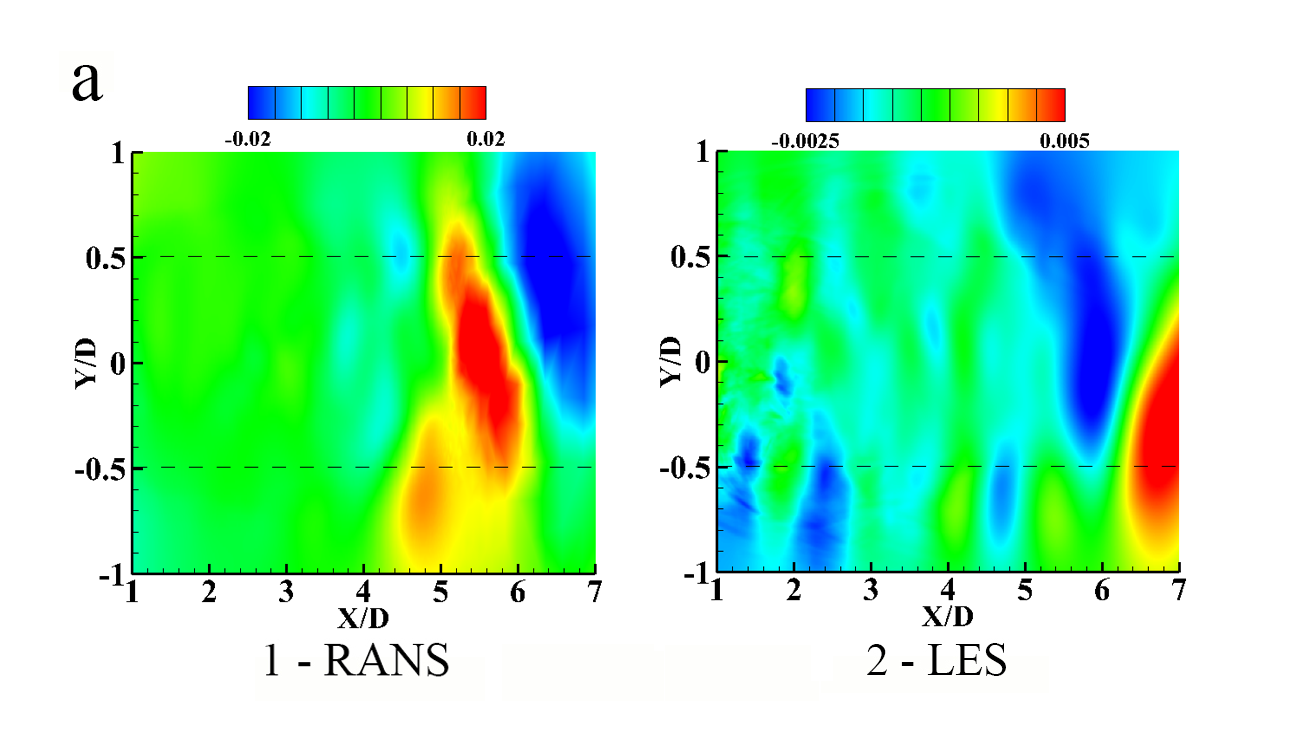}
	\includegraphics[scale=.7]{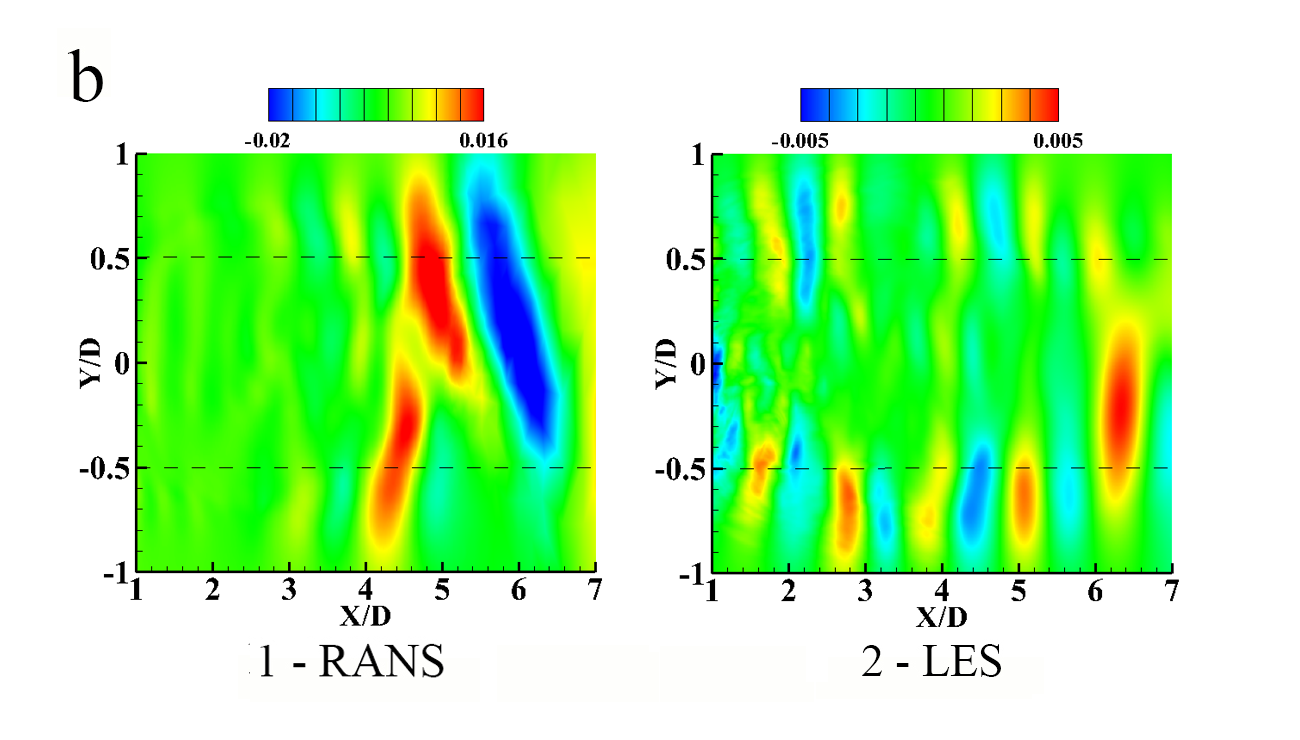}
	\includegraphics[scale=.7]{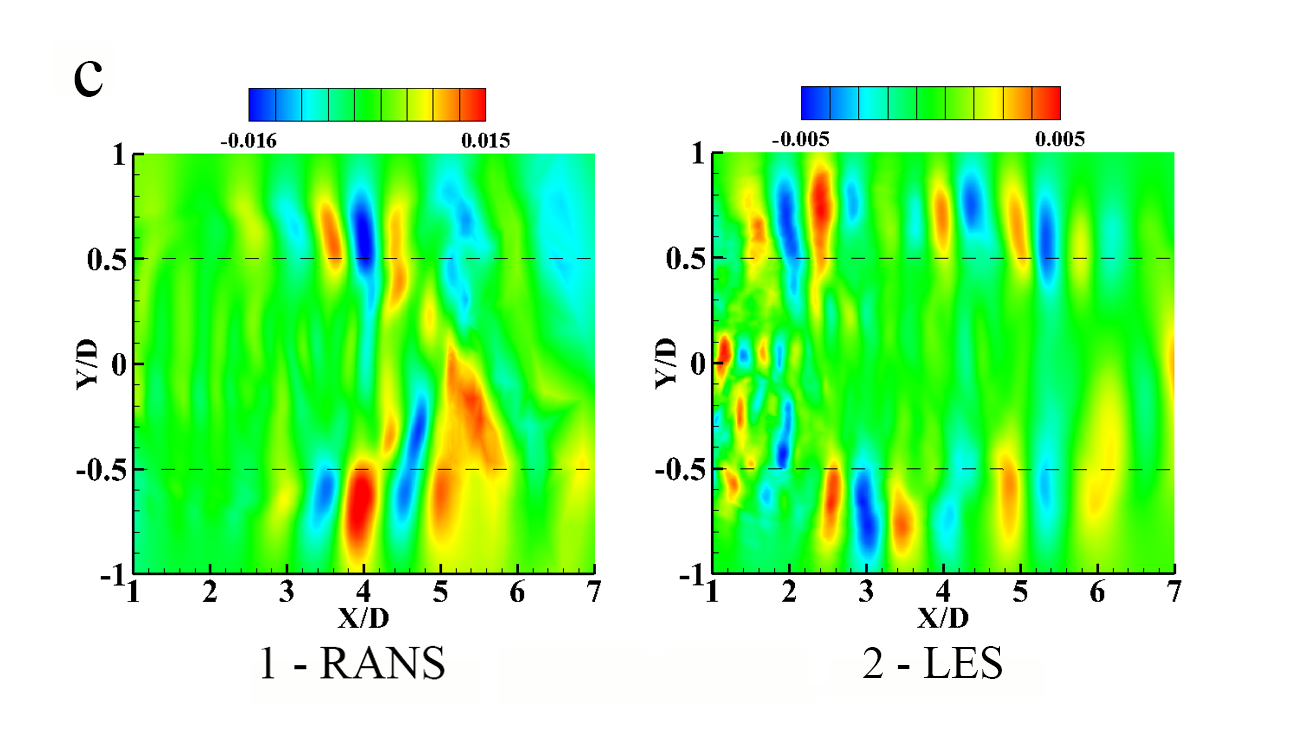}
	\includegraphics[scale=.7]{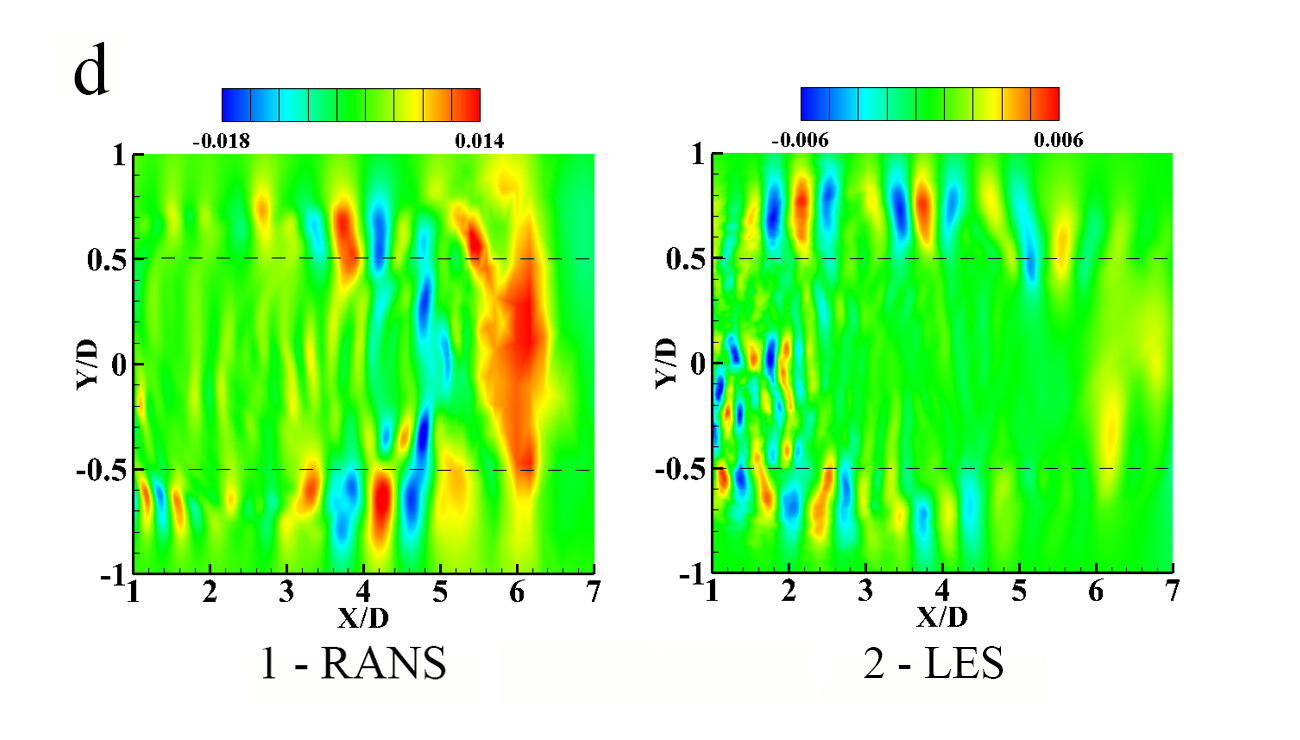}
	\includegraphics[scale=.7]{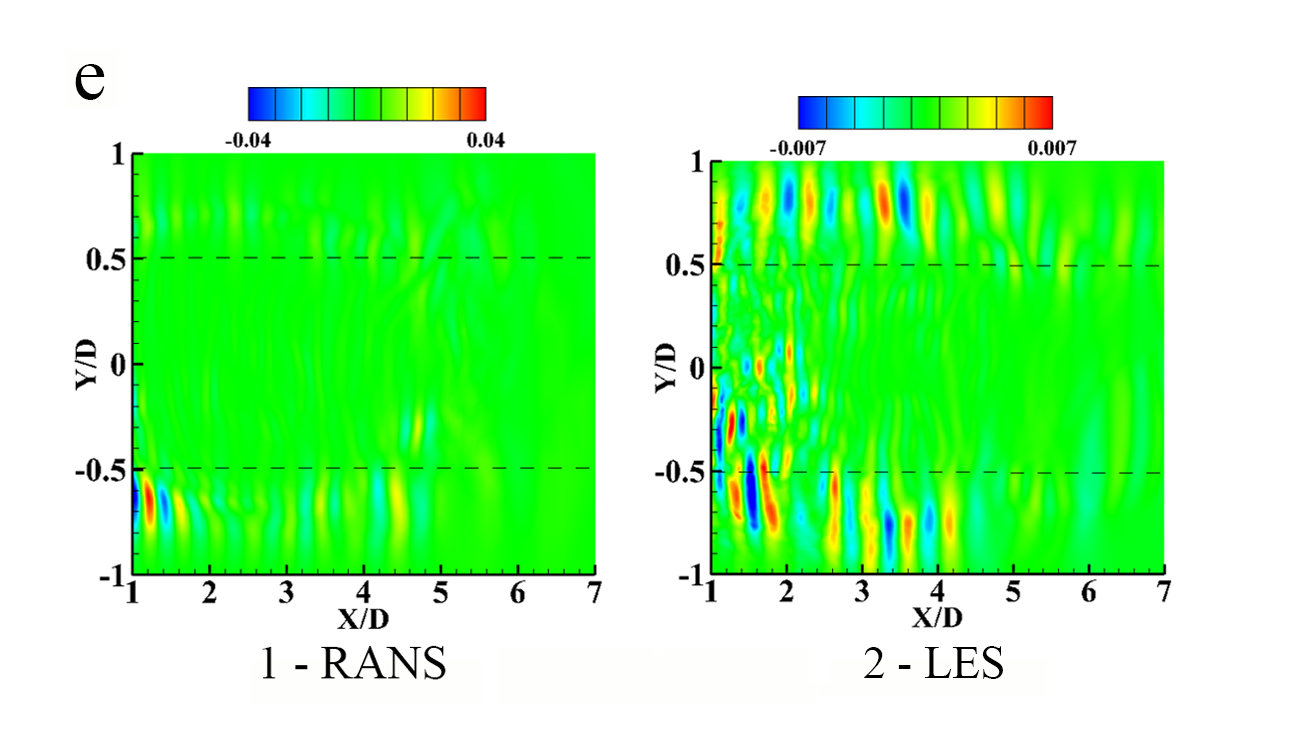}
	\caption{POD modes for transversal component of velocity at TSR2.4 (left RANS and right LES)}
	\label{FIG4_6}
\end{figure}

The spatial modes for the transversal component of velocity for TSR2.4 have been presented in Fig. \ref{FIG4_6}. The first transversal RANS mode, unlike that of the streamwise component of velocity, shows flow structures at the beginning of the wake region collapse point which is 5D behind the rotor and the initial part of the wake region is almost free of noticeable structures. Therefore, the wake region collapse leads to entrainment between the wake region and the far field area. It should be noted that due to the direction of transversal modes, the structures associated with these modes could be interpreted as the mixing process. Similar to POD modes of the streamwise component of velocity, with the increase in the mode number the decay of the flow structures is observed which is evident in Fig. \ref{FIG4_6}-b-1 and Fig. \ref{FIG4_6}-c-1 for the second and the third modes of transversal component of velocity respectively. In addition, the fourth and the fifth RANS modes presented in Fig. \ref{FIG4_6}-d-1 and \ref{FIG4_6}-e-1 respectively demonstrate the rolled-up symmetrical flow structures at the edges of the wake region identical to that of the streamwise component of velocity at the same mode numbers. 

Contours of spatial modes for the LES simulations reveals a quite similar pattern with that of the RANS. Similarly, as the LES mode number increases, the structures near the edge of the POD box appear more prominently, while the structures at the point of wake region collapse begin to disappear. Also, as the mode number increases from the first to the second mode at TSR2.4, it is observed that the mixing structures at point of wake region collapse begin to decay into structures of lower energy. In addition, as the mode number increases the vortical structures near the edges of the wake region rise into importance, while those at the point of wake region collapse become less noticeable which could be seen in Fig. \ref{FIG4_6}-c-2 to Fig. \ref{FIG4_6}-e-2. While the RANS results fail to demonstrate vortex dissipation at the wake region collapse area from third to forth modes (Fig. \ref{FIG4_6}-c-1 and Fig. \ref{FIG4_6}-d-1), LES successfully succeeds in capturing the physical expected behavior of structures breakup into smaller ones (Fig. \ref{FIG4_6}-c-2 and \ref{FIG4_6}-d-2).  It should be noted that this is in agreement with the observation for the modes of the streamwise component of velocity where first the structures due to the wake region collapse appear and gradually the rolled symmetrical structures start to appear. However, the main difference with the RANS modes, is that in the LES modes the structures due to the wake shed and entrainment ensuing from the rotor which are near the beginning of the wake region begin to rise into importance which are absent in the RANS modes. Moreover, contrary to the first mode of the streamwise component of velocity (Fig. \ref{FIG4_6}-a-2 and Fig. \ref{FIG4_6}-a-2), it is found that the first transversal mode manifests a convective behavior. This could be better understood by considering the time evolution of POD coefficients which will be discussed later on.   

\begin{figure*}
	\centering
	\includegraphics[scale=.3]{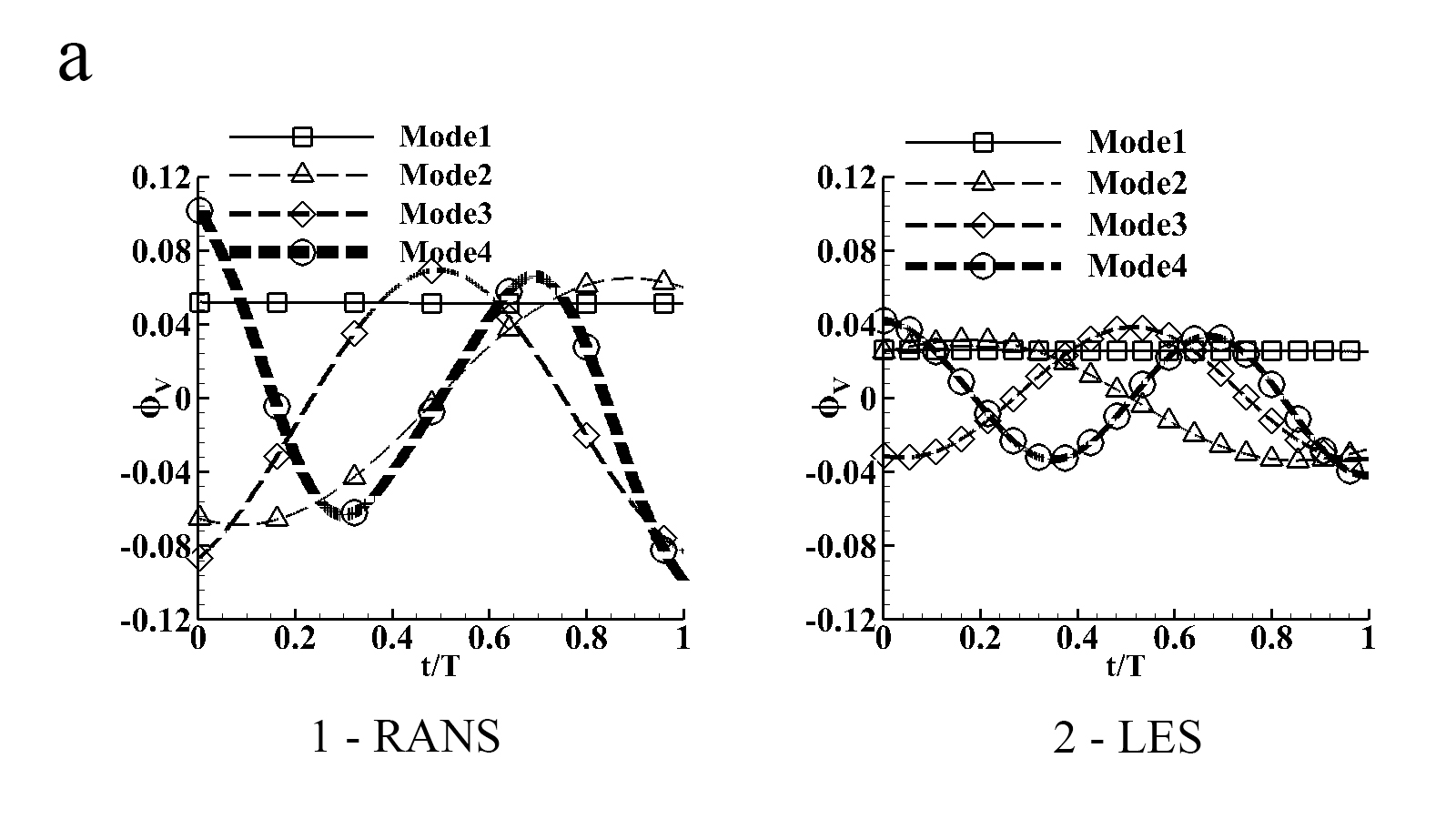}
	\includegraphics[scale=.3]{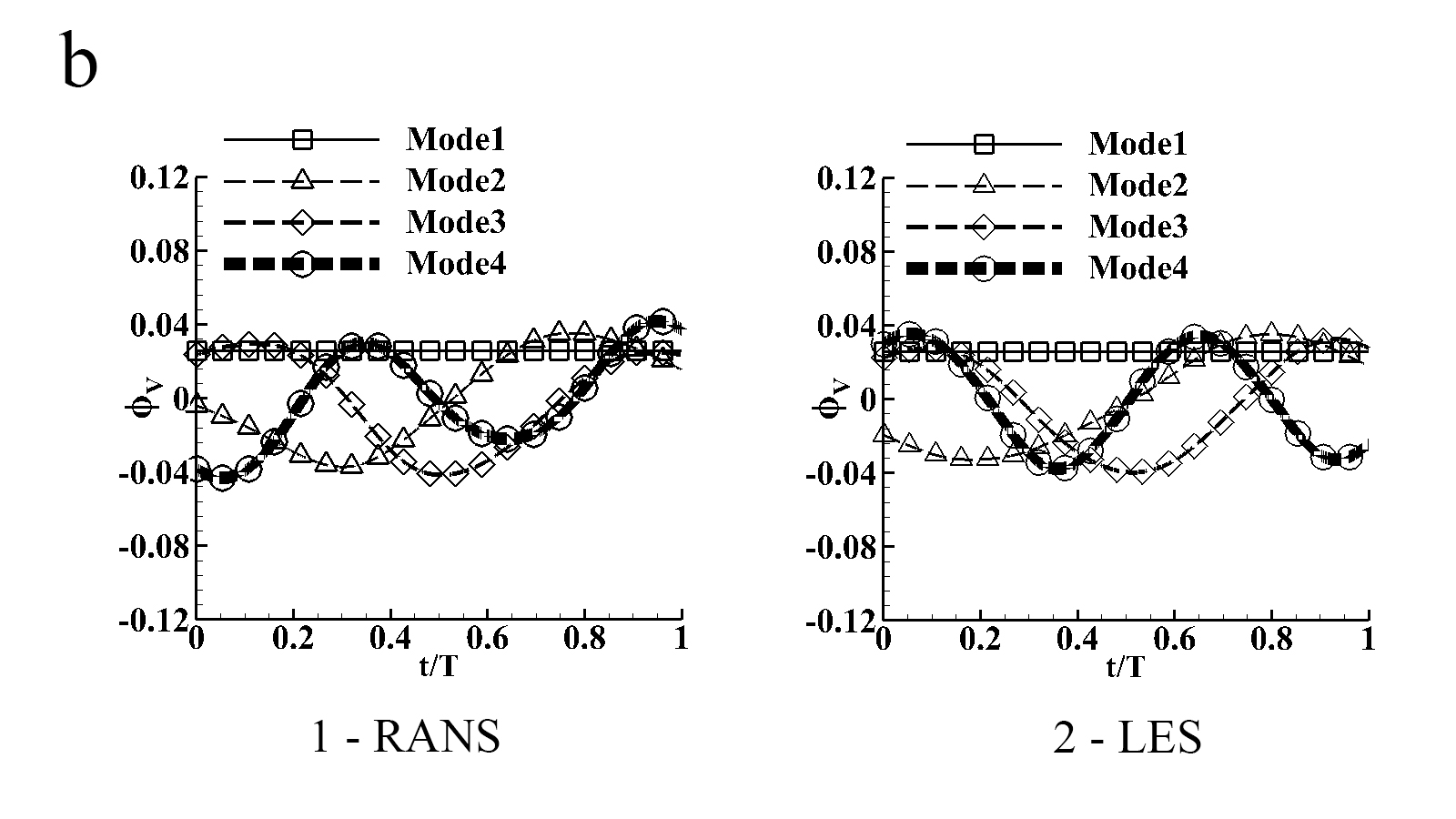}
	\includegraphics[scale=.3]{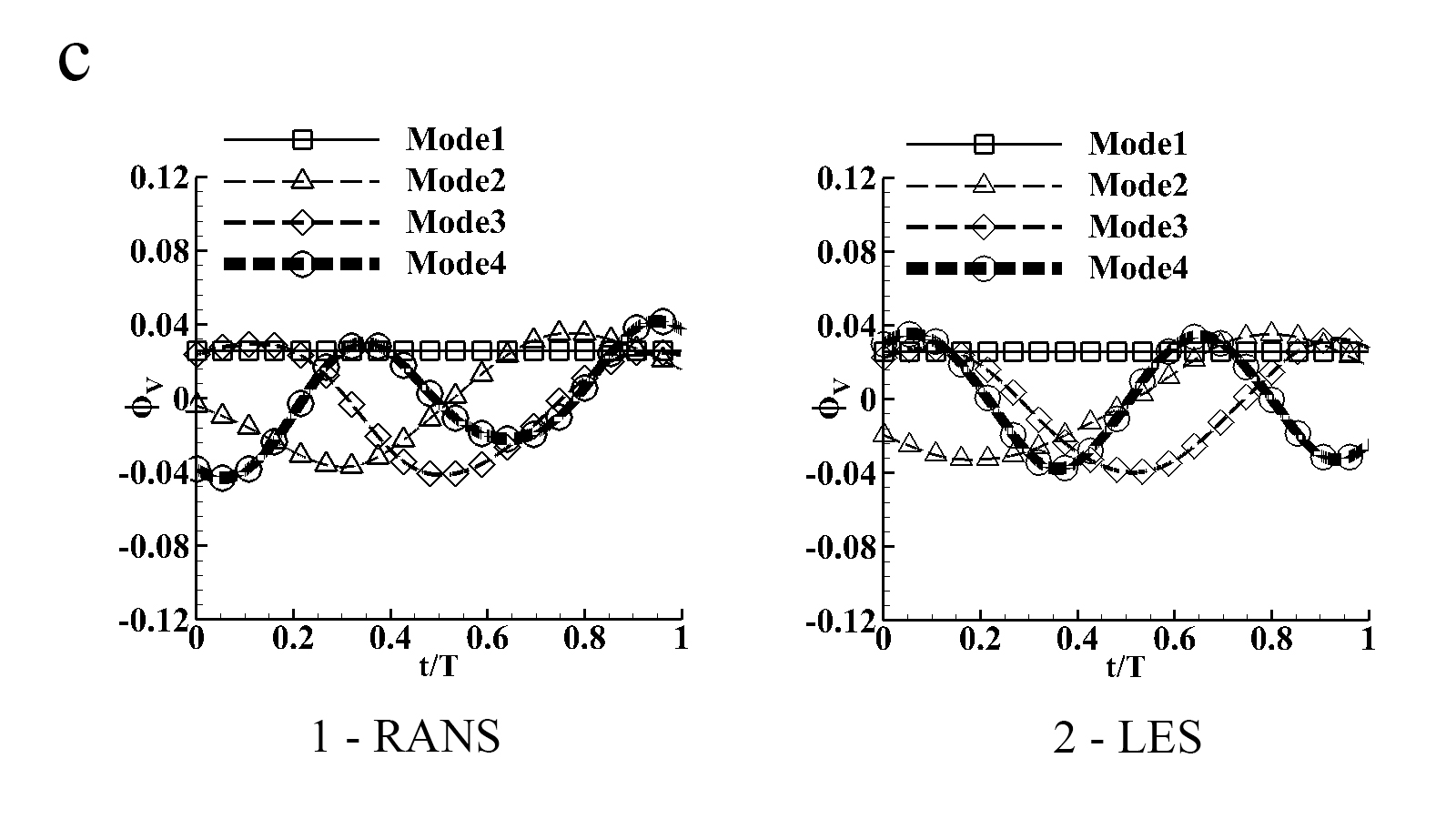}
	\caption{Time evolution of POD coefficents of the streamwise component of velocity for (a) TSR2.4 (b) TSR3.3 (c) TSR1.5}
	\label{FIG4_7}
\end{figure*}

\begin{figure*}
	\centering
	\includegraphics[scale=.3]{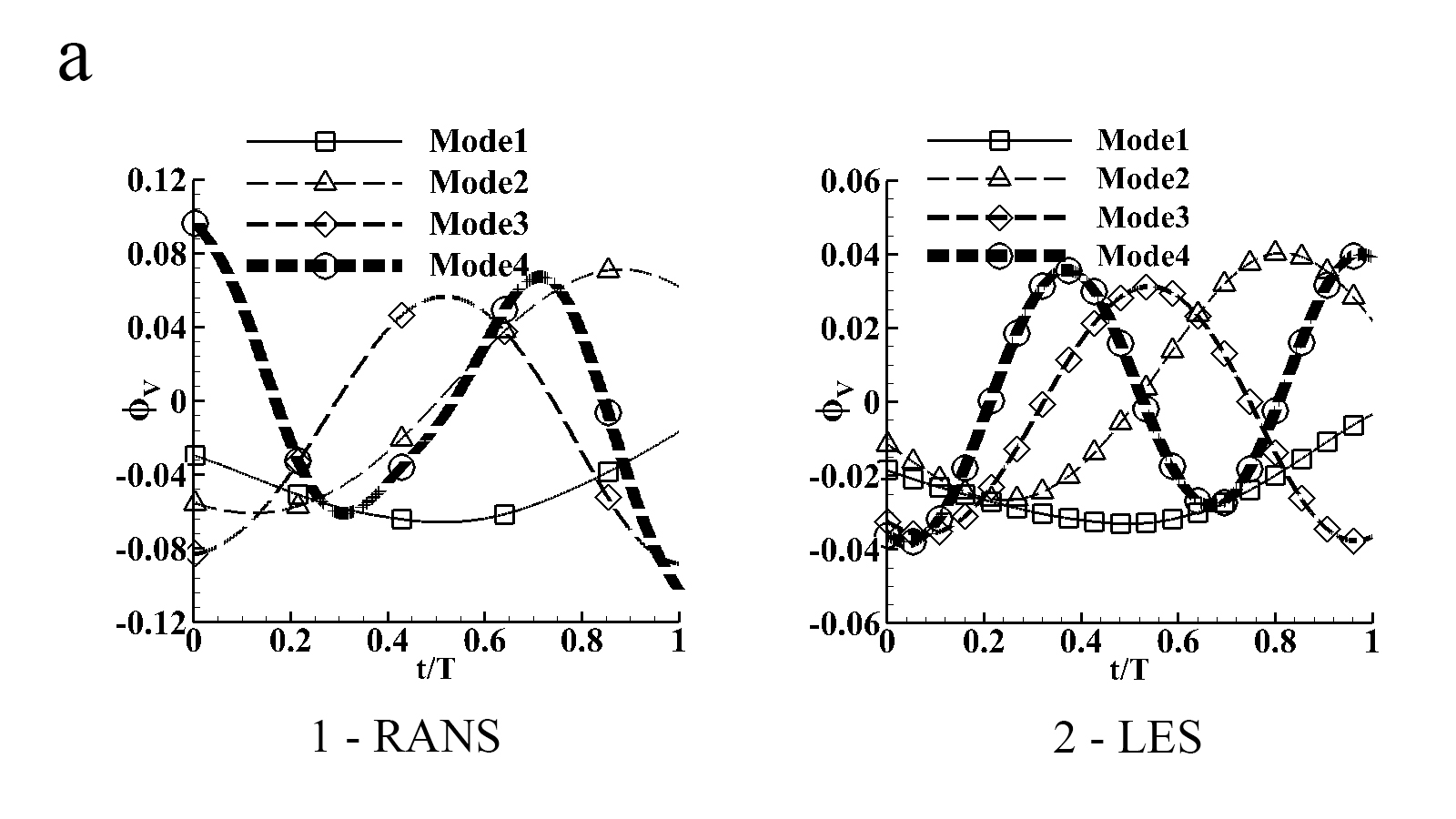}
	\includegraphics[scale=.3]{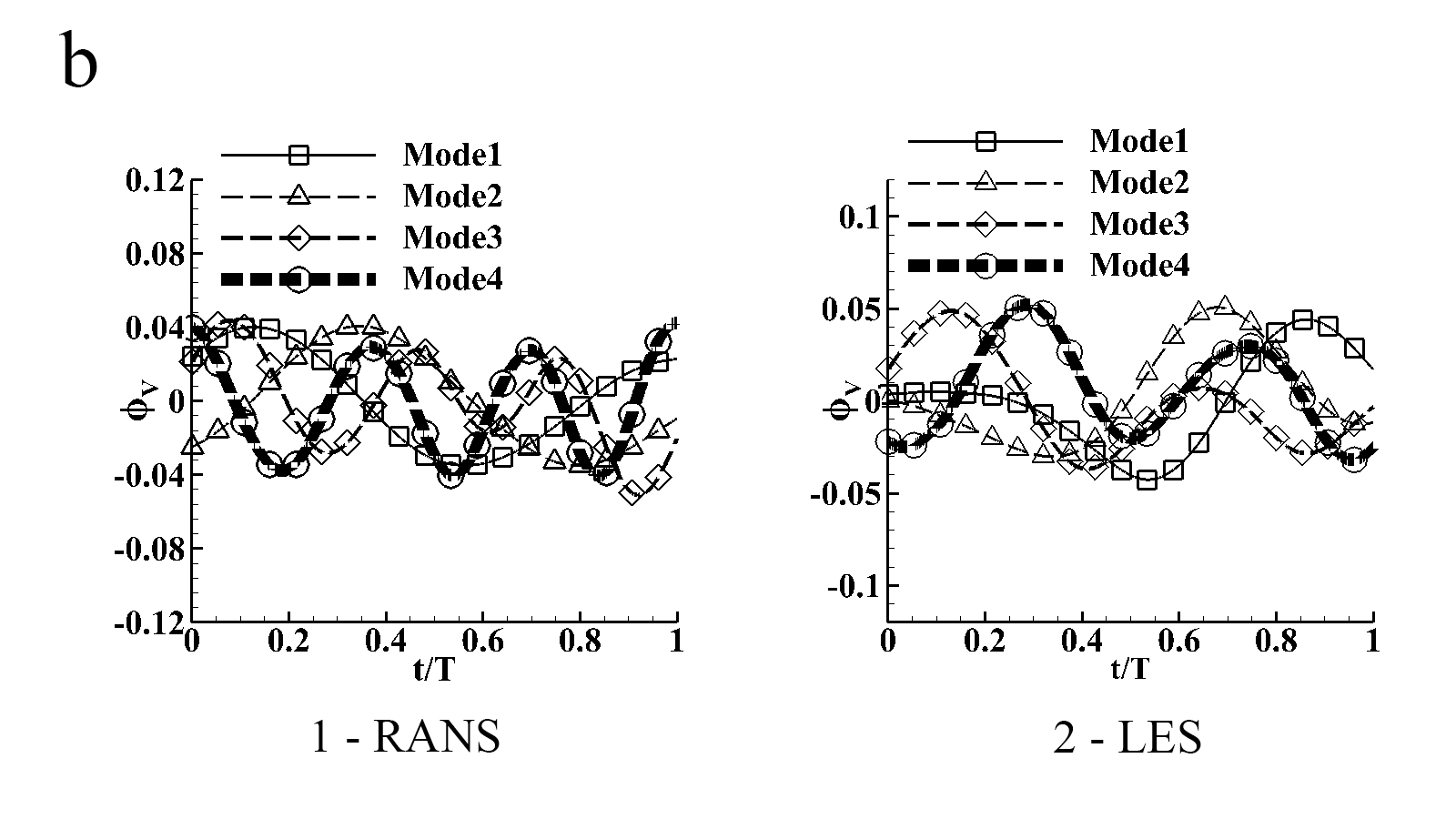}
	\includegraphics[scale=.3]{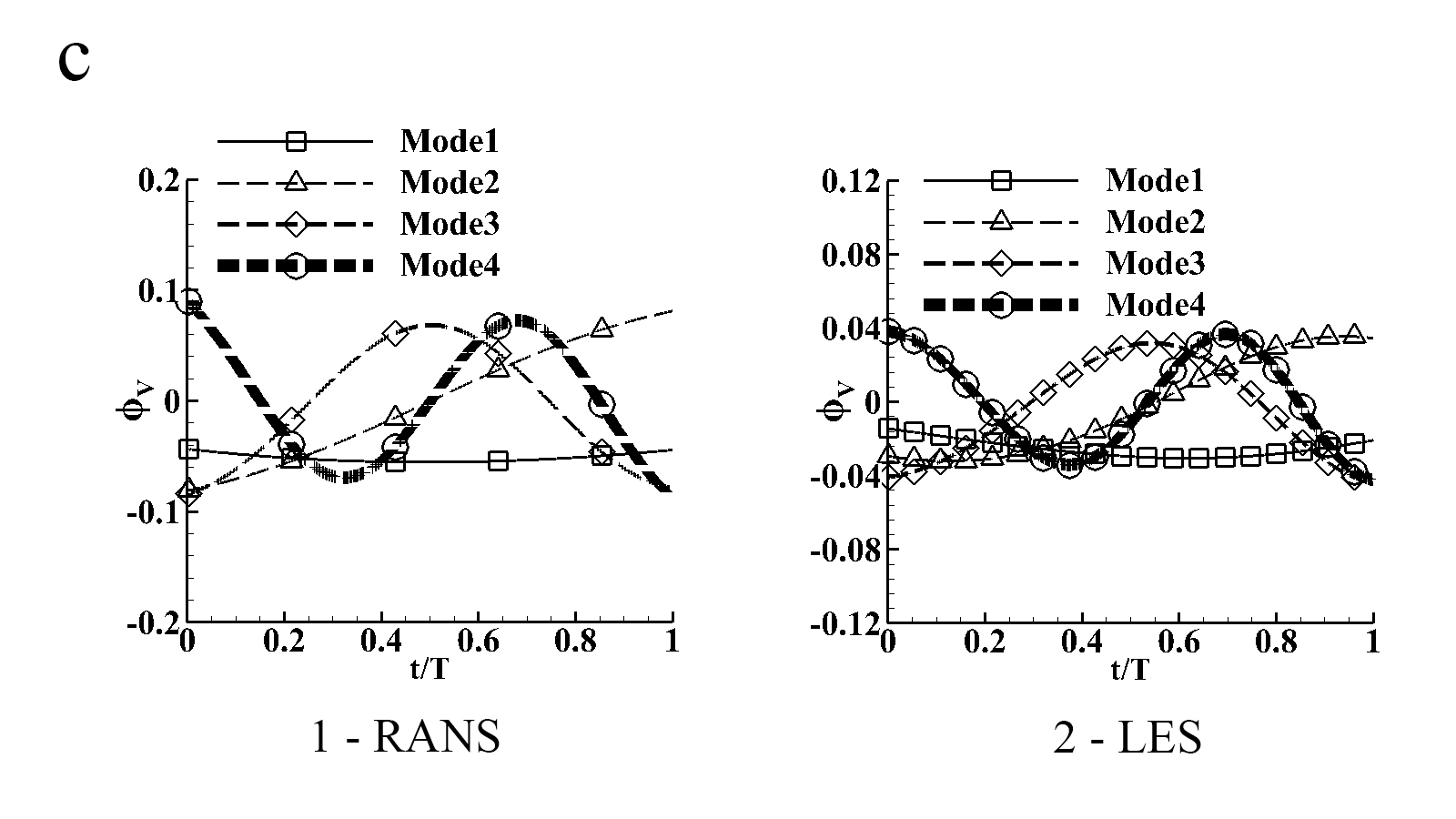}
	\caption{Time evolution of POD coefficents of the transversal component of velocity for (a) TSR2.4 (b) TSR3.3 (c) TSR1.5}
	\label{FIG4_8}
\end{figure*}

The time evolution of POD coefficients of the streamwise and transversal components of velocity have been presented in Fig. \ref{FIG4_7}. First, the RANS time evolution of POD coefficients will be explicated. As observed in Figs. \ref{FIG4_7}-a-1, \ref{FIG4_7}-b-1 and \ref{FIG4_7}-c-1 the first RANS time evolution of POD coefficient of the streamwise component of velocity at all TSRs demonstrate a constant trend proving a non-convective behavior of the first coefficient of time evolution of POD. Therefore, the trace of the first coefficient of time evolution of POD can be hardly observed in the consecutive modes. Following this, as the mode number increases, a rise in the mode frequency is observed. As a result, the convective behavior of the POD modes will be more significant with the increase in the mode numbers which is in congruent with the spatial modal behavior. Also, the frequency of the modes increases with a decrease in TSR as it is evident in Fig. \ref{FIG4_7}. Hence, as previously observed, low TSRs correspond to more oscillatory flow structures which are attributed to the larger blade loadings at higher wind velocities. Regarding the time evolution of POD coefficients of transversal component of velocity, it is observed in Fig. \ref{FIG4_8}-a that the first mode shows a transient trend especially, at higher TSRs. Therefore, contrary to the time evolution of POD coefficient for streamwise component of velocity, the first coefficient of time evolution of POD of the transversal component of velocity demonstrates a convective characteristic. Also, after the first coefficient of time evolution of POD, a striking resemblance between time evolution of POD coefficient of the both components of velocity is observed indicating that a large portion of modes energy has been almost equally distributed between the streamwise and the transversal components of velocity. In other words, similar flow structures exist in both U and V modes and the similarities become more considerable with the rise in the modes number. Also, the behavior of time evolution of POD coefficients at TSR3.3 and TSR2.4 reveals more resemblance with respect to TSR1.5 insinuating that the vortex shedding at the most loaded sections of the wind turbine, which is typical of low TSRs, leads to a change in the trend and an increase in the connectiveness of the POD modes because of the higher frequency of each mode compared to its counterpart at higher TSRs. Also, as seen in Fig. \ref{FIG4_8}-c, it should be noted that as the time evolution of POD coefficients of the transversal component of velocity demonstrate more evident similarity compared to that of the streamwise, it is concluded that the vortex shedding at lower TSRs affects the V modes more considerably. It is also worth noting that the time evolution of POD coefficients show a similar trend to that of flow past a stationary barrier as found in the work of Lengani et al. \citep{Lengani_2014}. Therefore, similar flow structures and modal behavior especially, at lower mode numbers could be expected for these two cases.

Generally, the LES and RANS time evolution of POD coefficients follow a similar trend. To illustrate, as anticipated from the spatial modes, the first time evolution of POD coefficients of the streamwise component of velocity at all TSRs reveal a non-convective behavior as seen in Fig. \ref{FIG4_7}-a-1, \ref{FIG4_7}-b-1 and \ref{FIG4_7}-c-1. In addition, with the increase in the mode number, the frequency of the time evolution of POD coefficients rises. However, at some modes some difference in phase could be observed. For instance, as seen in Fig. \ref{FIG4_7}-a, there is a $180^{\circ}$ phase difference between second LES and RANS time evolution of POD coefficient of the streamwise component of velocity at TSR2.4. Moreover, as shown in Fig. \ref{FIG4_7}-b, the RANS time evolution of POD coefficients TSR3.3 reveal the most similarity in trend to the LES results. This is attributed to the fact that at high TSRs the vortex shedding due to the lower blade loading leading to the more analogous modal behavior as also observed in spatial modes.

\section{Conclusion}
In this study the POD analysis of the wake shed behind a VAWT was performed by RANS and LES methods. To this end, WALE method, as suggested by Ma et al. \citep{Ma_2009}, was employed for SGS modelling. In addition, having examined various RANS models at the peak power, the well-known $k-\omega$ SST method was selected to perform the simulation for other TSRs. It should be noted that the LES resolution was found to be at least 82\% and the $y^{+}$ was well below unity. 

In the following step, the POD analysis was implemented on the wake behind the rotor. Firstly, the normalized cumulative modal energy was plotted for both LES and RANS simulations and it was shown that the RANS simulation results in higher modal energy at each mode number compared to that of the LES. Also, it was found that the RANS cumulative modal energy at all TSRs undergoes a considerably faster growth with respect to the LES energy modes. This was attributed to the capability of LES which resolves a significantly larger part of the flow domain leading to the presence of more flow structures than that of the RANS. Moreover, it was observed that the highest modal is directly related to the TSR. In other words, the highest TSR, for both the streamwise and transversal components of velocity, corresponds to the highest modal energy, while the lowest TSR leads to the lowest with respect to the other TSRs. 

In the second step, the spatial modes were presented and it was observed that the first spatial RANS mode for the streamwise component of velocity at all TSRs shows a resemblance to that of the LES as the first mode represents the area affected by the work extraction by the rotor. Also, with the increase in the mode number, more discrepancies among the TSR modes begin to appear which is in agreement with the trend shown in cumulative modal energy curves. Regarding the spatial modes for the transversal component of velocity, the contours revealed that both the RANS and LES first modes for TSRs 2.4 (rated) and 3.3 correspond to the flow structures due to the entrainment at the point of wake region collapse, while the first transversal and the successive modes at TSR1.5 demonstrates the vortex shedding due to the high blade loading at lower TSRs. As for the time evolution of POD coefficients, it was observed that both the RANS and LES modes demonstrate similar trends with the increase in the mode number. In this regard, it was concluded that the RANS method can capture the trend in the time evolution of POD coefficients, however, with different amounts.           

In conclusion, the RANS method can be employed for the purpose of the calculation of the wake region length and width or the extraction of the flow structures at the location of the wake region collapse. Nevertheless, a thorough POD analysis of the wake region cannot be performed by the RANS method as it is not capable of capturing the flow structures of smaller dimensions. As a result, for instance, as for an industrial or academic investigation of the wake region behind a VAWT, the RANS method can be employed to obtain an overall view of the wake region e.g. the wake region collapse length. In addition, the effect of the VAWT on the wake at a distance near the wake region collapse point can be accurately modelled by the RANS method. However, the internal area of the wake region especially, located near the rotor is associated with structures that RANS fails to capture indicating that at such points the employment of LES is more recommended if a deep understanding of the wake region is intended.

In summary, the main findings of the paper can be listed as follows:

\begin{itemize}
\item The RANS results demonstrate higher modal energy at each mode number compared to those of LES.
\item The first and second modes are well-captured by the RANS method owing to the fact the flow structures at these modes are mainly of the integral scale.
\item Moving on to the consecutive modes, the discrepancies start to rise into importance as the RANS accuracy ebbs for flow structures of smaller dimensions, which are typical of higher mode numbers due to the lower energy and higher frequency.
\item Time Evolution of the POD coefficients for RANS calculation at all TSRs demonstrate a very similar trend to that of the LES.
\end{itemize}

\section{Bibliography}

\bibliographystyle{cas-model2-names}

\bibliography{refs}


\end{document}